\documentclass[11pt,a4paper,twoside]{article}
\usepackage{a4}
\usepackage{color}
\usepackage{bm, amsmath, amssymb, amsthm} 

\usepackage{graphicx}
\usepackage{subfigure}

\def\calf{\mathcal{F}}
\def\<{\langle}
\def\>{\rangle}
\def\eps{\varepsilon}
\def\NN{\mathbb{N}}
\def\ZZ{\mathbb{Z}}
\def\RR{\mathbb{R}}
\def\tr{\operatorname{Tr\,}}
\def\id{\operatorname{id\,}}
\def\Div{\operatorname{div}}
\newcommand{\dtau}{\operatorname{d}\!\tau}

\def\eq{\hspace*{-1.5mm}&=&\hspace*{-1.5mm}}
\def\plus{\hspace*{-1.5mm}&+&\hspace*{-1.5mm}}

\def\dt{\partial_t}
\def\CC{{\mathbb C}}

\newtheorem{corollary}{Corollary}

\newtheorem{example}{Example}
\newtheorem{remark}{Remark}
\newtheorem{lemma}{Lemma}
\newtheorem{proposition}{Proposition}
\newtheorem{theorem}{Theorem}

\author{Vladimir Rovenski\thanks{
        E-mail: rovenski@math.haifa.ac.il
        }
        \ and \
        Leonid Zelenko\thanks{
        E-mail: zelenko@math.haifa.ac.il
        }\\
        {\small Mathematical Department, University of Haifa}
}

\title{The mixed Yamabe problem for harmonic foliations}

\date{}

\begin{document}

\maketitle

\begin{abstract}
The mixed scalar curvature of a foliated Riemannian manifold, i.e., an ave\-raged mixed sectional curvature,
has been considered by several geometers.
We explore the Yamabe type problem: to prescribe the constant mixed scalar curvature for a foliation
by a conformal change of the metric in normal directions only.
For a harmo\-nic foliation, we derive the leafwise elliptic equation and explore the corresponding nonlinear heat type equation.
We assu\-me that the leaves are compact submanifolds and the manifold is fibered instead of being foliated,
and use spectral parame\-ters of certain Schr\"{o}dinger operator to find solution, which is attractor of the equation.



\end{abstract}


\section*{Introduction}
\label{sec:main-res}

Geometrical problems of prescribing curvature-like invariants (e.g. the scalar curvature and the mean curvature)
of manifolds and foliations are popular for a long time, see \cite{cc2,rovwol,sw}.
There are many proofs of positive answer to the Yamabe problem:
gi\-ven a closed Rieman\-nian manifold $(M,g)$ of $\dim M\ge 3$,
find a metric conformal to $g$ with constant scalar curvature.
The study of this geometrical problem  was began by Yamabe (1960)
and completed by Trudinger, Aubin and Schoen (1986),
its solution is expressed in terms of the existence and multiplicity of solutions of a given elliptic PDE in the Riemannian manifold, see \cite{aub,r-o}.
Several authors developed an analog of the problem for CR-manifolds, see \cite{ho2012}, and its generalization
to contact (real or quaternionic) manifolds.
The problem when metrics of constant scalar curvature can be produced
on warped product manifolds has been studied in several articles, see~\cite{du}.

Let $(M,g)$ be endowed with a $p$-dimensional foliation~$\calf$.
Denote by ${\mathcal D}=T\calf$ the tangent distribution and
${\mathcal D}^\bot\,(\dim{\mathcal D}^\bot=n)$ the
orthogonal distribution (or the normal subbundle) of the tangent bundle $TM$.
In~\cite{bf}, a tensor calculus, adapted to the orthogonal splitting
\begin{equation}\label{E-TM-decomp}
 TM = {\mathcal D} + {\mathcal D}^\bot
\end{equation}
is developed to study the geometry of both the distributions and the ambient mani\-folds.
We~have $g = g_\calf + g^\perp$, where $g^\perp(X,Y):=g(X^\perp,Y^\perp)$
and $(\,\cdot\,)^\perp$ is the projection of $TM$ onto ${\mathcal D}^\bot$.
Obviously, biconformal metrics $\tilde g = v^{\,2} g_\calf + u^{\,2} g^\perp$\ $(u,v>0)$
preserve (\ref{E-TM-decomp}) and extend the class of conformal metrics (i.e., $u=v$).
Biconformal metrics (e.g. doubly-twisted products, introduced by Ponge and Reckziegel in \cite{pr})
have many applications in differential geometry, relativity, quantum-gravity, etc, see \cite{du}.
The~${\mathcal D}^\bot$- or ${\mathcal D}$-conformal metrics correspond to $v\equiv1$ or $u\equiv1$,
see \cite{rw-m}\,--\,\cite{rz2013}.

Using the natural representation of ${\rm O}(p)\times {\rm O}(n)$ on $TM$,
A.M. Naveira \cite{N1983} distinguished  thirty-six classes of Riemannian almost-product manifolds
$(M,g,{\cal D},{\cal D}^\bot)$; some of them are foliated, e.g.,
harmonic, totally geodesic, conformal, and Riemannian foliations.
Following this line of research, several geometers completed the geometric interpretation and gave examples for each class
of almost-product structures.
The simple examples of harmonic foliations are geodesic ones (e.g., parallel circles or winding lines on a flat torus,
and a Hopf field of great circles on the 3-sphere).
H.\,Rummler
characterized harmonic foliations by existence of an $\calf$-closed diffe\-rential $p$-form
that is transverse to $\calf$.
D.\,Sullivan's topological tautness condition is equivalent to the existence of a
metric on $M$ making a foliation harmonic, see~\cite{cc2}.

The components of the curvature of a foliation can be tangential, transversal, and mixed.
The tangential geometry of a foliation is infinitesimally
modeled by the tangent distribution to the leaves, while the transversal geometry
-- by the orthogonal distribution ${\mathcal D}^\bot$.
Prescribing the sign of tangential sca\-lar curvature has been studied for foliated spaces,
for example, there is no foliation of positive leafwise scalar curvature on any torus, see \cite{zh2013}.
The transversal scalar curvature is well studied for Riemannian foliations,
e.g. the ``transversal Yamabe problem", see~\cite{wz2013}.

The mixed scalar curvature, ${\rm S}_{\rm mix}$, for foliated (sub)manifolds
has been considered by several geometers, see \cite{bip,pp,wa1},
but its constancy (so called ``mixed Yamabe problem") is less stu\-died.
 In \cite{rz2012,rz2013}, we prescribed the sign of ${\rm S}_{\rm mix}$ using flows of ${\mathcal D}^\bot$-conformal metrics.
 In~this paper we explore the following Yamabe type problem:
\textit{Given a harmonic foliation $\calf$ of a Riemannian manifold $(M,g)$, find a ${\mathcal D}^\bot$-conformal
metric $\tilde g$ with leafwise constant~mixed scalar curvature.}
For a general foliation, the topology of the leaf through a point can change dramatically with the point;
this gives many difficulties in stu\-dying leafwise parabolic and elliptic equations.
 Therefore, in the paper (at least in the main results) we assume~that
\begin{equation}\label{E-2conditions}
 \calf\ \mbox{\rm is defined by an orientable fiber bundle}.
\end{equation}

The proof of main results is based on Sect.~1.2 (with variation formulae for geometrical quantities
under ${\mathcal D}^\bot$-conformal change of a metric),
Sect.~1.3 (with Proposition~\ref{T-main04} and Corollary~\ref{T-main0-Riemfol}),
Sect.~\ref{sec:evolPDEs} (about attractor of the nonlinear heat equation on a closed manifold
and about solution of its stationary equation with parameter)
and Sect.~\ref{sec:app} (about smooth dependence of a solution on a transversal parameter).

A slight change in the proof allows us to extend the main results for the case
when the prescribed mixed scalar curvature is not leafwise constant.

\section{Main results}
\label{subsec:prel}

The main results of the paper are the following.

\begin{theorem}\label{T-A}
Let $\calf$ be a harmonic and nowhere totally geodesic foliation
of a closed Riemannian manifold $(M, g)$ with condition~(\ref{E-2conditions}).
Then there exists a ${\mathcal D}^\bot$-conformal metric $\tilde g$ with leafwise
constant mixed scalar curvature.
\end{theorem}

If ${\mathcal D}^\bot$ is integrable than Corollary~\ref{T-main0-Riemfol} (see Sect.~\ref{sec:mainAB}) is applicable.
In particular case of codimen\-sion-one foliations, we have the following.

\begin{corollary}\label{C-A}
Let $\calf$ be a codimension-one harmonic and nowhere totally geodesic foliation
of a Riemannian manifold $(M, g)$ with condition~(\ref{E-2conditions}).
Then there exists a ${\mathcal D}^\bot$-conformal metric $\tilde g$ with
leafwise constant Ricci curvature in the normal direction.
\end{corollary}

There are examples of foliations of codimension $>1$ with minimal,
not totally geode\-sic leaves on (compact) Lie groups with left-invariant
metrics, see \cite{ty84}; further, the metric can be chosen to be bundle-like with respect to $\calf$.
Such foliations have leafwise constant mixed scalar curvature.

\begin{theorem}\label{T-B}
Let $\calf\ (\dim\calf>1)$ be a totally geodesic foliation of a closed Riemannian manifold $(M, g)$
with condition~(\ref{E-2conditions}) and integrable normal distribution.
Then there exists a ${\mathcal D}^\bot$-conformal metric $\tilde g$ with leafwise constant  mixed scalar curvature.
\end{theorem}

\textbf{1.1. Preliminaries}.
Denote by $R(X,Y)=\nabla_Y\nabla_X-\nabla_X\nabla_Y+\nabla_{[X,Y]}$ the~curvature tensor of Levi-Civita connection.
The sec\-tional curvature $K(X, Y)=g(R(X,Y)X,Y)$, where $X\in T\calf,\,Y\in{\mathcal D}^\bot$ are unit vectors,
is called \textit{mixed}.
It regulates (through the Jacobi equation) the deviation of leaves along the leaf geodesics.
Foliations with constant mixed sectional curvature play an important role in differential geo\-metry,
but are far from being classified.
Examples are $k$-nullity foliations on Riemannian manifolds which are totally geodesic,
relative nullity foliations, which determine a ruled structure of submanifolds in space forms,
foliations produced by Reeb vector field on Sasakian manifolds, etc.
 Totally geodesic foliations on complete manifolds with $K_{\rm mix}\!\equiv0$~split.
For a $k$-dimensional totally geodesic foliation with $K_{\rm mix}\equiv1$ on a closed manifold $M^{n+k}$,
we have the Ferus's inequality $k<\rho(n)$, where $\rho(n)-1$ is the number of
linear independent vector fields on a sphere $S^{n-1}$, see~\cite{rov-m}.

The~\textit{mixed scalar curvature} is an averaged mixed sectional curvature,
\[
 {\rm S}_{\rm mix} =\sum\nolimits_{j=1}^n\sum\nolimits_{a=1}^p K({\mathcal E}_j, {E}_a),
\]
and is independent of the choice of a local orthonormal frame $\{{\mathcal E}_j,\,{E}_a\}_{j\le n,\,a\le p}$ of $TM$
adapted to ${\mathcal D}^\bot$ and $ T\calf$, see \cite{r2010,rov-m,wa1}.
If either ${\mathcal D}^\bot$ or $ T\calf$ is one-dimensio\-nal and tangent to a unit vector field $N$,
then ${\rm S}_{\rm mix}$ is the Ricci curvature in the $N$-direction.

Let ${\mathfrak X}_M$ be the module over $C^\infty(M)$ of all vector fields on $M$,
and ${\mathfrak X}^\bot$ and ${\mathfrak X}^\top$ the modules of all vector fields on ${\mathcal D}^\bot$ and $ T\calf$, respectively. The extrinsic geometry of a foliation is related to the second fundamental form of the leaves,
 $h(X,Y) = (\nabla_X Y)^\perp$, where $X,Y\in{\mathfrak X}^\top$,
and its invariants (e.g., the mean curvature $H=\tr_g\,h$).
Special classes of foliations such as totally geodesic, $h = 0$ (with the simplest extrinsic geometry);
totally umbilical, $h=(H/p)\,g_\calf$; and harmonic, $H=0$,
have been studied by many geometers, see survey in \cite{rov-m}.
Let $h^\bot$ be the second fundamental form of ${\mathcal D}^\bot$,
$H^\bot=\tr_g\,h^\bot$ the mean curvature,
and $T$ the integrability tensor of ${\mathcal D}^\bot$. We~have
\begin{equation}\label{E-def-B}
 2\,h^\bot(X,Y) = (\nabla_X\,Y+\nabla_Y\,X)^\top,\quad
 2\,T(X,Y)=[X,\,Y]^\top,\quad X,Y\in {\mathfrak X}^\bot.
\end{equation}
The formula in \cite{wa1}, for foliations reads as
\begin{equation}\label{E-PW}
 {\rm S}_{\,\rm mix}=\|H^\bot\|^2-\|h^\bot\|^2+\|T\|^2+\|H\|^2-\|h\|^2 +\Div(H^\bot + H).
\end{equation}
 We calculate norms of tensors using local adapted basis~as
\begin{eqnarray*}
 \|h^\bot\|^2 =\sum\limits_{i,j}\|h^\bot({\mathcal E}_i,{\mathcal E}_j)\|^2,\ \
 \|h\|^2 =\sum\limits_{a,b}\|h({E}_a,{E}_b)\|^2,\ \
 \|T\|^2 =\sum\limits_{i,j}\|T({\mathcal E}_i,{\mathcal E}_j)\|^2\,.
\end{eqnarray*}

\begin{example}\rm({Constant mixed scalar curvature on doubly-twisted products}).
The \textit{doubly twisted product} of Riemannian manifolds $(B,g_\calf)$ and $(F, g^\perp)$,
is a manifold $M=B\times F$ with the metric $g = v^2\,g_\calf + u^2\,g^\perp$,
where $v,u\in C^\infty(B\times F)$ are positive functions.
It is called the \textit{doubly warped product} of $(B,g_\calf)$ and $(F, g^\perp)$
if the warping functions $v$ and $u$ only depend on the points of $B$ and $F$, respectively.

The \textit{leaves} $B\times\{y\}$ of a doub\-ly-twisted product $B\times_{(v,\,u)} F$
and the \textit{fibers} $\{x\}\times F$ are totally umbilical. We have
\[
 h = -(\nabla^\perp \log v)\,g_\calf,\quad h^\bot = -(\nabla^\top\log u)\,g^\perp.
\]
By the above, $H=-n\,\nabla^\perp\log v$, $H^\bot=-p\,\nabla^\top\log u$, and
\[
 \|H\|^2-\|h\|^2=(n^2-n)\|\nabla^\perp v\|^2/v^2,\quad
 \|H^\bot\|^2-\|h^\bot\|^2=(p^2-p)\|\nabla^\top u\|^2/u^2.
\]
Next we derive
\begin{eqnarray*}
 \Div\,H \eq -p\,(\Delta^\top\,u)/u -(p^2-p)\|\nabla^\top u\|^2/u^2,\\
 \Div\,H^\bot \eq -n\,(\Delta^\perp v)/v -(n^2-n)\|\nabla^\perp v\|^2/v^2,
\end{eqnarray*}
where $\Delta^\top$ is the leafwise Laplacian and $\Delta^\perp$ is the fiberwise Laplacian.
Substituting in (\ref{E-PW}) with $T=0$, we obtain the formula
\[
 {\rm S}_{\,\rm mix} = -n\,(\Delta^\top u)/u -p\,(\Delta^\perp\,v)/v\,.
\]
Let $B$ be a closed manifold.
Given a positive function $v\in C^\infty(B\times F)$, define a leafwise Schr\"{o}dinger operator
$\mathcal{H}=-\Delta^\top-\beta$, where $\beta = \frac pn\,(\Delta^\perp\,v)/v$.
For any compact leaf, the spectrum of $\mathcal{H}$ is discrete,
the least eigenvalue $\lambda_0$ is isolated from other eigenvalues,
and the eigenfunction $e_0$ (called the ground state) can be chosen positive, see Sect.~\ref{sec:evolPDEs}.
Since $\mathcal{H}(e_0)=\lambda_0\,e_0$, a doubly-twisted product $B\times_{(v,\,e_0)} F$
has leafwise constant mixed scalar curvature equal to $n\lambda_0$.
\end{example}

By Lemma~\ref{L-btAt}, ${\mathcal D}^\bot$-conformal changes of the metric preserve harmonic foliations.

We focus on the mixed Yamabe problem for harmonic foliations,
which amounts to finding a positive solution of the leafwise elliptic equation,
see Proposition~\ref{P-Yam1},
\begin{equation}\label{E-Yam1-init}
 -n\,(\Delta^\top u +\beta^\top\,u) = -2\,H^\bot(u) +\widetilde{\rm S}_{\rm mix}\,u +\|h\|^2_{g}\,u^{-1} -\|T\|^2_{g}\,u^{-3}\,,
\end{equation}
where
 $\beta^\top = \frac1n(\,\|T\|^2_{g}-\|h\|^2_{g} - {\rm S}_{\rm mix})$\,,
and a leafwise constant $\widetilde{\rm S}_{\rm mix}$ corresponds
to a ${\mathcal D}^\bot$-confor\-mal metric $\tilde g$.
Proposition~\ref{T-02} allows us to reduce (\ref{E-Yam1-init}) to the case of $H^\bot=0$.

\begin{example}\label{Ex-main0-folB}\rm
The global structure of totally geodesic foliations with integrable normal bundle
(i.e., ${\mathcal D}^\bot$ is tangent to a foliation $\calf^\bot$) has been studied in \cite{bh1983}:
the universal cover $\tilde M$ is topologically a product $\tilde F\times\tilde F^\bot$
of universal covers of the leaves of both foliations, $\calf$ and $\calf^\bot$.
Let $\calf$ be a totally geodesic foliation with integrable normal bundle
of a closed Riemannian manifold $(M, g)$ with conditions~(\ref{E-2conditions}) and $H^\bot=0$.
Then $\Psi_1=\Psi_2=0$, and (\ref{E-Yam1-init}) becomes the linear elliptic equation on $F$,
\begin{equation}\label{E-Yam1-init-Linear}
 -\Delta^\top u -(\beta^\top+\Phi)\,u  = 0\,,
\end{equation}
where $\beta^\top=-\frac 1n\,{\rm S}_{\rm mix}$.
Suppose that ${\rm S}_{\rm mix}\ne{\rm const}$ and $\Phi={\rm const}$.
Then ${\cal H}\,(u_*)=\Phi\,u_*$, where $u_*=e_0$ and $\Phi=\lambda_0$
for the Schr\"{o}dinger operator ${\cal H}=-\Delta^\top -\beta^\top$.
Assuming $\nabla^\bot u_{\,|F}=0$, continue $u_*$ smoothly on $M$.
Thus, the~mixed scalar curvature of the Riemannian manifold
$(M,\tilde g=g^\top+u_*^2 g^\perp)$ is~$n\Phi$.
\end{example}

\begin{proposition}\label{T-02}
Let $\calf$ be a foliation of a closed Riemannian manifold $(M, g)$ with condition~(\ref{E-2conditions}).
Then there exists a smooth function $u>0$ on $M$ such that $H^\bot=0$
for the metric $\tilde g=g_\calf + u^2 g^\perp$.
\end{proposition}

\noindent\textbf{Proof}.
Recall the equality for any $X,Y\in{\mathfrak X}^\bot$ and $U,V\in{\mathfrak X}^\top$, see~\cite{rov-m},
\begin{eqnarray}\label{E-genricA}
\nonumber
 g(R(U,X)V,Y) \eq g(((\nabla_U\,C)_V-C_VC_U)(X),Y)\\
 \plus g(((\nabla_X\,A^\top)_Y - A^\top_X A^\top_Y)(U),V),
\end{eqnarray}
where the co-nullity operator $C:  T\calf\times TM\to {\mathcal D}^\bot$
is defined~by
 $C_U(X)=-(\nabla_{\!X} U)^\perp$\ $(U\in{\mathfrak X}^\top,\ X\in{\mathfrak X}_M)$.
Note~that
\begin{eqnarray*}
 \sum\nolimits_j\,g((\nabla_U\,C)_V({\mathcal E}_j),\,{\mathcal E}_j)
 \eq\sum\nolimits_j\nabla_U (g(C_V({\mathcal E}_j),\,{\mathcal E}_j))\\
 \eq\nabla_U \big(g\,\big(\sum\nolimits_j\,h({\mathcal E}_j,{\mathcal E}_j),\,V\big)\big) = g(\nabla_U H^\bot,\,V).
\end{eqnarray*}
Thus, tracing (\ref{E-genricA}) over $\mathcal D$ and taking the antisymmetric part, we obtain
$d^\top H^\bot = 0$, where the $2$-form $d^\top H^\bot$ is defined by
\[
 2\,d^\top H^\bot(U,V) = g(\nabla_U H^\bot, V)-g(\nabla_V H^\bot, U)\quad (U,V\in{\mathfrak X}^\top).
\]
Then we apply Lemma~\ref{L-rw} given below.
\qed

\begin{lemma}[\rm see Theorem~1.1 in \cite{rovwol}\,]\label{L-rw}
Let $\calf$ be a foliation of a closed Riemannian manifold $(M, g)$ with condition~(\ref{E-2conditions}),
and $d^\top H^\bot=0$. Then the Cauchy's problem
\[
 \dt g = -(2/p)\,(\Div^\top H^\bot)\,g^\perp,\quad g_0=g,
\]
has a unique solution $g_t\ (t\ge0)$ that converges as $t\to\infty$ to a metric with $H^\bot=0$.
\end{lemma}

\textbf{1.2. ${\mathcal D}^\bot$-conformal change of a metric}.
 We shall find how various geometrical quantities are transformed under
${\mathcal D}^\bot$-confor\-mal change of a metric.
 The~Weingarten operator $A^\bot_U$ of ${\mathcal D}^\bot$ and the
skew-symmetric operator $T^\sharp_U$, where $U\in{\mathfrak X}^\top$, are given~by
\[
 g(A^\bot_U(X),\,Y) = g(h^\bot(X,Y),\ U),\quad
 g(T^\sharp_U(X),\,Y) = g(T(X,Y),\ U).
\]

\begin{lemma}
\label{L-btAt}
Given a foliation $\calf$ on $(M,g = g_\calf + g^\perp)$, and $\phi\in C^1(M)$, define a new metric
$\tilde g=g_\calf + e^{\,2\/\phi} g^\perp$. Then
\begin{eqnarray}
\label{E-Ak-tau}
 \tilde h\,^\top \eq e^{-2\/\phi} h,\quad \tilde H\,^\top = e^{-2\/\phi} H, \\
\label{E-Ak-U}
 \tilde h\,^\bot \eq e^{\,2\/\phi}\big(h^\bot -(\nabla^\top\phi)\,g^\perp\big),\quad
 \tilde H\,^\bot = H^\bot -n\,\nabla^\top\phi, \\
\label{E-Ak-T}
 \tilde A\,^\bot_U \eq A^\bot_U -U(\phi)\,\id^\perp,\quad
 \tilde T_U^\sharp = e^{-2\/\phi}\,T_U^\sharp\quad(U\in{\mathfrak X}^\top).
\end{eqnarray}
Hence, ${\mathcal D}^\bot$-conformal variations preserve total umbilicity, harmonicity, and total geo\-desy of~$\,\calf$,
and preserve total umbilicity of the normal distribution $\,{\mathcal D}^\bot$.
\end{lemma}

\noindent\textbf{Proof}.
 The \textit{Levi-Civita connection} $\nabla$ of a metric $g$ is given by the known formula
\begin{eqnarray}\label{eqlevicivita}
\nonumber
 &&\quad 2\,g(\nabla_X\,Y, Z) = X g(Y,Z) + Y g(X,Z) - Z g(X,Y) \\
 && +\,g([X, Y], Z) - g([X, Z], Y) - g([Y, Z], X)\quad (X,Y,Z\in{\mathfrak X}_M).
\end{eqnarray}
Formula (\ref{E-Ak-tau})$_1$ follows from (\ref{eqlevicivita}):
\begin{eqnarray*}
 2e^{2\/\phi} g(\tilde\nabla_U V, X) \eq 2\tilde g(\tilde\nabla_U V, X)\\
 \eq -X g(U,V) -g([U,X],V) -g([V,X],U) =2g(\nabla_U V, X).
\end{eqnarray*}
We deduce (\ref{E-Ak-tau})$_2$ using
 $\tilde H\,^\top =e^{-2\/\phi}\sum\nolimits_a h({E}_a,{E}_a)=e^{-2\/\phi}H$.
From $\tilde T=T$ and
\begin{eqnarray*}
 g(\tilde T_U^\sharp(X), Y)\eq e^{-2\/\phi}\tilde g(\tilde T_U^\sharp(X), Y)
 =e^{-2\/\phi}\tilde g(T(X,Y), U)\\
 \eq e^{-2\/\phi} g(T(X,Y), U) =e^{-2\/\phi}g(T_U^\sharp(X), Y)
\end{eqnarray*}
formula (\ref{E-Ak-T})$_2$ follows.
By (\ref{eqlevicivita}), for any $X,Y\in{\mathfrak X}^\bot$ and $U\in{\mathfrak X}^\top$ we have
\begin{equation}\label{E-compute1}
 g(\widetilde\nabla_X Y,\,U) = e^{\,2\/\phi}\,g(\nabla_X Y,\,U) -e^{\,2\/\phi}\,U(\phi)\,g(X,Y)
 -(e^{\,2\/\phi}-1)\,g(T(X,Y),\,U).
\end{equation}
From this, skew-symmetry of $T$ and (\ref{E-def-B}), we deduce (\ref{E-Ak-U})$_{1}$.
Then we get (\ref{E-Ak-T})$_{1}$ using
\begin{eqnarray*}
 e^{\,2\/\phi}g(\tilde A\,^\bot_U(X),Y) \eq \tilde g(\tilde A\,^\bot_U(X),Y) = \tilde g(\tilde h\,^\bot(X,Y), U)\\
  \eq e^{\,2\/\phi}\big(g(A^\bot_U(X),Y) - U(\phi)\,g(X,Y)\big).
\end{eqnarray*}
Similarly, we prove (\ref{E-Ak-T})$_{2}$:
\[
 e^{\,2\/\phi}g(\tilde T^\sharp_U(X),Y) =\tilde g(\tilde T^\sharp_U(X),Y)
 =\tilde g(\tilde T(X,Y), U) = g(T^\sharp_U(X),Y).
\]
The orthonormal bases of ${\mathcal D}^\bot$ in both metrics are related by $\tilde {\mathcal E}_j = e^{-\phi}{\mathcal E}_j$.
To show this we calculate for any $j\le n$,
\begin{equation}\label{E-epsconfframe}
 1=\tilde g(\tilde {\mathcal E}_j, \tilde {\mathcal E}_j) =e^{\,2\/\phi}\, g(e^{-\phi} {\mathcal E}_j, e^{-\phi} {\mathcal E}_j) = g({\mathcal E}_j, {\mathcal E}_j).
\end{equation}
By (\ref{E-Ak-U})$_1$, we have
\[
 \tilde h\,^\bot(\tilde {\mathcal E}_j,\tilde {\mathcal E}_j)=e^{-2\phi}\,\tilde h\,^\bot({\mathcal E}_j, {\mathcal E}_j)= h^\bot({\mathcal E}_j, {\mathcal E}_j)-(\nabla^\top\phi)\,g({\mathcal E}_j, {\mathcal E}_j).
\]
From this and the definition $H^\bot=\tr_{g}\,h^\bot$, the equality (\ref{E-Ak-U})$_2$ follows.
\qed

\begin{remark}\label{R-btAt}\rm By Lemma~\ref{L-btAt},
for a leafwise constant $\phi$ we obtain $\tilde h\,^\bot = e^{\,2\/\phi}h^\bot$ and $\tilde H\,^\bot = H^\bot$.
Hence, ${\mathcal D}^\bot$-scalings of \,$g$ preserve harmonicity and total geodesy of~$\,{\mathcal D}^\bot$.
\end{remark}

\begin{proposition}\label{P-Yam1}
The mixed scalar curvature of a harmonic foliation $\calf$ under ${\mathcal D}^\bot$-conformal change of a~metric
$\tilde g = g_\calf + u^{2}g^\perp$, where $u>0$ is a smooth function, is transformed according to the formula
\begin{equation}\label{E-Yam1-ini}
  ({\rm S}_{\,\rm mix} -\widetilde{\rm S}_{\rm mix})\,u
  =  n\,\Delta^\top\,u -2\,H^\bot(u) +\|h\|^2_{g}\,(u^{-1}-u) -\|T\|^2_{g}\,(u^{-3}-u).
\end{equation}
If, in particular, $\,u$ is leafwise constant $($i.e., $\tilde g$ is a ${\mathcal D}^\bot$-scaling of $g)$, then we have
\[
 \widetilde{\rm S}_{\rm mix} = {\rm S}_{\rm mix} -(u^{-2} -1)\|h\|^2_{g} +(u^{-4} - 1)\|T\|^2_{g}\,.
\]
\end{proposition}

\noindent\textbf{Proof}. By Lemma~\ref{L-btAt}, we have
\begin{eqnarray}\label{E-CC-norms1}
\nonumber
 \|\tilde h^\top\|^2_{\tilde g} \eq e^{-2\phi}\|h\|^2_{g},\ \
 \|\tilde T\|^2_{\tilde g} = e^{\,-4\phi}\|T\|^2_{g},\ \
 \|\tilde h^\bot\|^2_{\tilde g} = \|h^\bot\|^2_{g} +n\|\nabla^\top\phi\,\|^2_{g} -2\,H^\bot(\phi),\\
 \|\tilde H\,^\bot\|^2_{\tilde g} \eq \|H^\bot\|^2_{g} +n^2\|\nabla^\top\phi\,\|^2_{g} -2\,n\,H^\bot(\phi),\quad
 \widetilde\Div\,^\top\tilde H\,^\bot = \Div^\top H^\bot -n\,\Delta^\top\phi\,.
\end{eqnarray}
Indeed, the formulae for $\|\tilde h\,^\top\|^2_{\tilde g}$ and $\|\tilde T\|^2_{\tilde g}$ follow from
\begin{eqnarray*}
 \|\tilde h\,^\top\|^2_{\tilde g} \eq\sum\nolimits_{a,b,i} \tilde g(\tilde h\,^\top(E_a, E_b),\tilde{\mathcal E}_i)^2
 = e^{4\,\phi}\sum\nolimits_{a,b,i} g(e^{-2\,\phi} h(E_a, E_b), e^{-\phi}{\mathcal E}_i)^2\\
 \eq e^{-2\,\phi}\sum\nolimits_{a,b,i} g(h(E_a, E_b),{\mathcal E}_i)^2 =e^{-2\,\phi}\|h\|^2_{g},\\
\|\tilde T\|^2_{\tilde g}\eq\sum\nolimits_{a,b,i}\tilde g(\tilde T(\tilde{\mathcal E}_i,\tilde{\mathcal E}_j),E_a)^2
 =\sum\nolimits_{a,b,i} g(T(e^{-\phi}\tilde{\mathcal E}_i, e^{-\phi}\tilde{\mathcal E}_j), E_a)^2\\
 \eq e^{-4\,\phi}\sum\nolimits_{a,b,i} g(T({\mathcal E}_i, {\mathcal E}_j),E_a)^2
 =e^{-4\,\phi}\|T\|^2_{g}\,.
\end{eqnarray*}
Formula for $\widetilde\Div\,^\top\tilde H\,^\bot$ follows from
$g(\tilde\nabla_a\,U,E_a)=g(\nabla_a\,U,E_a)$ for $U\in{\mathfrak X}^\top$ and
\[
 \widetilde\Div\,^\top\tilde H\,^\bot = \sum\nolimits_{a}\tilde g(\tilde\nabla_a \tilde H\,^\bot, E_a)
 \Div^\top H^\bot -n\,\Div^\top(\nabla^\top\phi)\,.
\]
From
\begin{eqnarray*}
 \|\tilde h\,^\bot\|^2_{\tilde g}\eq\sum\nolimits_{a,i,j}\tilde g(\tilde h\,^\bot(\tilde{\mathcal E}_i,\tilde{\mathcal E}_i),E_a)^2
 =\sum\nolimits_{a,i,j}
 \big(g(h^\bot({\mathcal E}_i,{\mathcal E}_i)-(\nabla^\top\phi)\,g({\mathcal E}_i,{\mathcal E}_i),E_a)\big)^2\\
 \eq \|h^\bot\|^2_{g} -2\,g(H^\bot,\nabla^\top\phi) +n\,\|\nabla^\top\phi\,\|^2_g,\\
 \|\tilde H\,^\bot\|^2_{\tilde g} \eq g(\tilde H\,^\bot,\tilde H\,^\bot)
 = g(H^\bot -n\nabla^\top\phi\,,\, H^\bot -n\nabla^\top\phi)\\
 \eq \|H^\bot\|^2_g -2\,n\,g(H^\bot,\nabla^\top\phi) +n^2\|\nabla^\top\phi\,\|^2_g\,
\end{eqnarray*}
the formulae for $\|\tilde h\,^\bot\|^2_{\tilde g}$ and $\|\tilde H\,^\bot\|^2_{\tilde g}$ follow.
Then, using (\ref{E-CC-norms1}), $\tilde{\mathcal E}_i=e^{\,-\phi}{\mathcal E}_i$, and
\[
 \widetilde{\rm S}_{\rm mix} = \sum\nolimits_{i}\tilde r(\tilde{\mathcal E}_i,\,\tilde{\mathcal E}_i)
 =e^{\,-2\/\phi} \sum\nolimits_{i}\tilde r({\mathcal E}_i,{\mathcal E}_i),
\]
we obtain the formula
\begin{equation}\label{E-Kmixphi-h}
 \widetilde{\rm S}_{\rm mix} = {\rm S}_{\rm mix} {-}n\big(\Delta^\top\phi {+}\|\nabla^\top\phi\|^2_g\big)
 {+}2 H^\bot(\phi) {+}(e^{-4\phi}{-}1)\|T\|^2_{g} -(e^{-2\phi}{-}1)\|h\|^2_{g}\,.
\end{equation}
Substituting $\phi=\log u$ and
 $\nabla^\top\phi = u^{-1}\nabla^\top u$,
 $\Delta^\top\phi = u^{-1}\Delta^\top\,u -u^{-2}\,\|\nabla^\top u\|^2_g$
into (\ref{E-Kmixphi-h}) yields the required formula (\ref{E-Yam1-ini}), which is equivalent to (\ref{E-Yam1-init}).
\qed

\smallskip

\textbf{1.3. Proof of the main results}.\label{sec:mainAB}
 Proposition~\ref{T-02} allows us to assume $H^\bot=0$.
 Then we associate with (\ref{E-Yam1-init}) the leafwise parabolic equation
 with a leafwise constant $\widetilde{\rm S}_{\rm mix}$
\begin{equation}\label{E-evol-u1}
 \dt u -\Delta^\top u -(\beta^\top +\widetilde{\rm S}_{\rm mix}/n)\,u = (\|h\|^2/n)\,u^{-1} -(\|T\|^2/n)\,u^{-3}\,.
\end{equation}
We shall study asymptotic behavior of solutions to (\ref{E-evol-u1}) with appropriate initial data
using spectral parameters of a leafwise Schr\"{o}\-dinger operator
\[
 \mathcal{H}^\top = -\Delta^\top-\beta^\top\,.
\]
The least eigenvalue $\lambda^\top_0$ of $\mathcal{H}^\top$ is simple and obeys the inequalities
\[
 \lambda^\top_0\in[-\max\nolimits_{\,F}\beta^\top, \  -\min\nolimits_{\,F}\beta^\top]\,,
\]
its eigenfunction $e_0$ (called the ground state) may be chosen positive, see Sect.~\ref{sec:evolPDEs}.
By~(\ref{E-2conditions}) and results in Sect.~\ref{sec:app}, the leafwise constant $\lambda^\top_0$ and $e_0$ are smooth on~$M$.

 Assume $h\ne0$ and $\Phi<n\,\lambda_0^\top$ and consider the functions (compare with Sect.~\ref{sec:evolPDEs}),
\begin{eqnarray}\label{E-biquad}
\nonumber
 \phi^\top_-(y)= -(n\lambda_0^\top-\Phi) y +\min\nolimits_{\,F}\,(\|h\|^2 e_0^{-2})y^{-1}
 -\max\nolimits_{\,F}\,(\|T\|^2 e_0^{-4}) y^{-3},\\
 \phi^\top_+(y)= -(n\lambda_0^\top-\Phi) y +\max\nolimits_{\,F}\,(\|h\|^2 e_0^{-2})y^{-1}
 -\min\nolimits_{\,F}\,(\|T\|^2 e_0^{-4}) y^{-3}.
\end{eqnarray}
If the discriminant
 $D={\min\nolimits_{\,F}\,(\|h\|^4 e_0^{-4})}-4\,(n\,\lambda^\top_0 -\Phi)\,\max\limits_{F}(\|T\|^2 e_0^{-4})>0$,
each of (\ref{E-biquad}) has four real roots (two of them are positive). Their maximal (positive)~roots
\begin{eqnarray*}
 y^\top_- = \Big(\frac{\min\limits_{F}(\|h\|^2 e_0^{-2}) +({\min\limits_{F}(\|h\|^4 e_0^{-4})
 -4(n\lambda_0^\top-\Phi)\max\limits_{F}(\|T\|^2 e_0^{-4})})^{1/2}}{2\,(n\lambda_0^\top-\Phi)}\Big)^{\frac12},\\
 y^\top_+ = \Big(\frac{\max\limits_{F}(\|h\|^2 e_0^{-2}) +({\max\limits_{F}(\|h\|^4 e_0^{-4})
 -4(n\lambda_0^\top-\Phi)\min\limits_{F}(\|T\|^2 e_0^{-4})})^{1/2}}{2(n\lambda_0^\top-\Phi)}\Big)^{\frac12},
\end{eqnarray*}
obey the inequalities
$y^\top_- < y_3^\top < y^\top_+$,
where $y_3^\top$ is the maximal root of $(\phi^\top_-)'(y)$,
\[
 y_3^\top=\Big(\,
 \frac{({\min\limits_{\,F}(\|h\|^4 e_0^{-4})+12(n\lambda_0^\top-\Phi)\max\limits_{\,F}(\|T\|^2 e_0^{-4})})^{1/2}
 -\min\limits_{\,F}(\|h\|^2 e_0^{-2})}{2(n\lambda_0^\top-\Phi)}\Big)^{\frac12}\,.
\]
For a positive function $f\in C(F)$ define
$\delta_f:=(\min\nolimits_{\,F}f)/(\max\nolimits_{\,F}f)\in(0,\,1]$.

\begin{proposition}\label{T-main04}
Let $\calf$ be a harmonic and nowhere totally geodesic foliation
on a Riemannian manifold $(M, g)$ with condition~(\ref{E-2conditions}) and $H^\bot=0$.
Then for any leafwise constant $\Phi\in C^\infty(M)$ obeying the inequalities (along any leaf $F$)
\begin{equation}\label{E-lambdaPhi}
 n\,\lambda^\top_0 -\delta^{-4}_{e_0}\,\frac{\min_{\,F}\,\|h\|^4}{4\,\max_{\,F}\,\|T\|^2}
 <\Phi<n\,\lambda^\top_0\,,
\end{equation}
there exists a unique $u_*$ in the set $\{\tilde u\in C^\infty(M):\ \tilde u > e_0\,y^\top_3\}$
such that the mixed scalar curvature of $\tilde g = g_\calf + u_*^2\,g^\perp$ is $\Phi$.
 Moreover,
 $y^\top_-\le u_*/e_0\le y^\top_+$
and $u_*=\lim\limits_{\,t\to\infty} u(\cdot\,,t)$, where $u$ solves (\ref{E-evol-u1}) with
$\widetilde{\rm S}_{\rm mix}=\Phi$, does not depend on the value $u(\cdot\,,0)=u_0> e_0\,y^\top_3$.
\end{proposition}

\noindent\textbf{Proof}.
By Theorem~\ref{contdeponqeigfunct} (in Sect.~\ref{sec:app}),
the leafwise constant $\lambda_0^\top$ (the least eigenfunction of $\mathcal{H}^\top$)
and its leafwise eigenvector $e_0$ are smooth, i.e., they belong to $C^\infty(M)$.
If $M$ is closed there then exist many $\Phi$'s obeying (\ref{E-lambdaPhi}), e.g.
$n\,\lambda^\top_0-\eps$ for small enough~$\eps>0$.

By conditions, any leaf $F_0$ has an open neighborhood
diffeomorphic to $F\times\RR^n$ and $F_0=F\times\{0\}$.
Since $F_q=F\times\{q\}$ are compact minimal submanifolds,
their volume form $d\,{\rm vol}_F=|g_{|F}|^{1/2}{\rm d} x$
does not depend on $q\in\RR^n$, see~\cite{rov-m}.
Thus,~vector bundles $\{L_2(F_q)\}_{q\in\RR^n}$ and $\{H^k(F_q)\}_{q\in\RR^n}$
coincide with products $L_2\times\RR^n$ and $H^k\times\RR^n$.

Let $\Phi$ obey (\ref{E-lambdaPhi}) and let $u_0>e_0\,y^\top_3$ hold. We shall use the notation
\[
 \beta=\beta^\top+\Phi/n,\quad \lambda_0=\lambda^\top_0-\Phi/n,\quad
 \Psi_1 = \|h\|^2/n,\quad \Psi_2=\|T\|^2/n.
\]
Then
(\ref{E-evol-u1}) with $\widetilde{\rm S}_{\rm mix}=\Phi$ becomes (\ref{Cauchy}),
while (\ref{E-lambdaPhi}) follows from
\[
 n\,\lambda^\top_0 - (1/4)\,{\min_{\,F}\,(\|h\|^4 e_0^{-4})}/{\max_{\,F}\,(\|T\|^2 e_0^{-4})}
 <\Phi<n\,\lambda^\top_0\,,
\]
which becomes (\ref{E-lambda0-cond1}).
Hence, the claim follows from Theorem~\ref{thattract1} (in Sect.~\ref{sec:evolPDEs}).
\qed

\begin{corollary}\label{T-main0-Riemfol}
Let $\calf$ be a harmonic and nowhere totally geodesic foliation of a Riemannian manifold $(M, g)$ with condition~(\ref{E-2conditions}),
integrable normal subbundle ${\mathcal D}^\bot$ and $H^\bot=0$.
Then for any leafwise constant $\Phi\in C^\infty(M)$ such that $\Phi<n\lambda^\top_0$
there exists a~unique positive function $u_*\in C^\infty(M)$ such that (along any leaf $F$)
\[
 ({n\,\lambda_0^\top-\Phi})^{-1}{\min_{\,F}\,(\|h\|^2 e_0^{-2})} \le u_*/e_0
 \le ({n\,\lambda_0^\top-\Phi})^{-1}{\max_{\,F}\,(\|h\|^2 e_0^{-2})}\,,
\]
and the mixed scalar curvature of the metric $\tilde g = g_\calf + u_*^2\,g^\perp$ is $\Phi$.
\end{corollary}

\noindent\textbf{Proof}. This is similar to the proof of Proposition~\ref{T-main04}.
Since $\Psi_2\equiv0$ and $\lambda_0>0$, each of
 $\phi_-=-\lambda_0+\Psi_1^- y^{-1}$
 and
 $\phi_+=-\lambda_0+\Psi_1^+ y^{-1}$
has one positive root $y_1^-=(\Psi_1^-/\lambda_0)^{1/2}$ and $y_1^+=(\Psi_1^+/\lambda_0)^{1/2}$,
see also Example~\ref{Ex-2const}(c).
\qed

\smallskip
\noindent\textbf{Proof of Theorem~\ref{T-A}}.
 By Proposition~\ref{T-02}, there exists a metric $g_1$, ${\mathcal D}^\bot$-conformal to $g$, for which $H^\bot=0$
(the mean curvature of ${\mathcal D}^\bot$). By Lemma~\ref{L-btAt}, the equality $H=0$ is preserved for $g_1$.
 By Proposition~\ref{T-main04}, there exists a metric $\tilde g$, ${\mathcal D}^\bot$-conformal to $g_1$, for which
$\widetilde{\rm S}_{\rm mix}$ is leafwise constant; moreover, $H=0$ holds.
\qed

\smallskip
\noindent\textbf{Proof of Theorem~\ref{T-B}}.
 By Corollary~\ref{T-main0-Riemfol}, there is a metric $g_1$ that is  ${\mathcal D}^\bot$-conformal to $g$, for which $H^\bot=0$.
 By Lemma~\ref{L-btAt}, $h=0$ is preserved for $g_1$.
 Since $T=0$, equation (\ref{E-Yam1-init}) reads as the eigenproblem
 $\mathcal{H}(u) = \frac1n\,\widetilde{\rm S}_{\rm mix}\,u$,
where
$\mathcal{H}=-\Delta^\top-\beta$ is a leafwise Schr\"{o}\-dinger operator
on $(M,g_1)$ with potential $\beta = -\frac1n\,{\rm S}_{\rm mix}(g_1)$.
Let $e_0>0$ be the ground state of $\mathcal{H}$ with the least eigenvalue $\lambda^\top_0$
(leafwise constant).
Thus, the metric $\tilde g = g_\calf + e_0^2\,g_1^\perp$
 has $\widetilde{\rm S}_{\rm mix}=n\lambda^\top_0$;
moreover, the equality $h=0$ is preserved for~$\tilde g$.
\qed

\section{Results for the nonlinear heat equation}
\label{sec:evolPDEs}

Let $(F, g)$ be a smooth closed $p$-dimensional Riemannian manifold (e.g., a~leaf of a compact foliation)
with the Riemannian distance $d(x,y)$.
Functional spaces over $F$ will be denoted without writing $(F)$, e.g., $L_2$ instead of $L_2(F)$.
 Let $H^k$ be the Hilbert space of Sobolev real functions of order $k$ on $F$
with the inner product $(\,\cdot,\cdot\,)_{k}$ and the norm $\|\cdot\|_k$.
In~particular, $H^0=L_2$ with the product $(\,\cdot,\cdot\,)_{0}$ and the norm~$\|\cdot\|_0$.
 Denote by $\|\cdot\|_{C^k}$ the norm in the Banach space $C^k$ for $1\le k<\infty$,
and $\|\cdot\|_C$ for $k=0$.
In~local coordinates $(x_1,\dots, x_p)$ on~$F$, we have
$\|f\|_{C^k}=\max_{F}\max_{|m|\le k}|d^m f|$, where $m\ge0$ is the multi-index of order
$|m|=\sum_{i}m_i$ and $d^m$ is the partial derivative (in fact, a finite atlas of $F$ must be considered).

\begin{proposition}[{\rm Scalar maximum principle, see \cite[Theorem~4.4]{ck1}}]\label{P-weak-max}
Let $X_t$ and $g_t$ be smooth families of vector fields and metrics
on a closed Riemannian manifold $F$, and $f\in C^\infty(\RR\times[0,T))$.
Suppose that $u:F\times[0,T)\to\RR$ is a $C^\infty$ solution to
\[
 \dt u \ge\Delta_t u - X_t(u) +f(u,t),
\]
and let $y:[0,T]\to\RR$ solve the Cauchy's problem for ODEs: $y\,' = f(y(t),t)$,\ $y(0)=y_0$.
If~$u(\cdot\,, 0)\ge y_0$, then $u(\cdot\,,t)\ge y(t)$ for all $t\in[0, T)$.
\end{proposition}

\textbf{2.1. The nonlinear heat equation}.
We are looking for stable solutions of the elliptic equation, see (\ref{E-Yam1-init}) with $H^\bot=0$,
\begin{equation}\label{Cauchy-stat1}
 -\Delta\,u -\beta\,u = \Psi_1(x)\,u^{-1} -\Psi_2(x)\,u^{-3}\,,
\end{equation}
where $\Psi_1>0$, $\Psi_2\ge0$ and $\beta$ are smooth functions on $F$.
To study (\ref{Cauchy-stat1}), we shall look for attractors of the Cauchy's problem for the {nonlinear heat equation},
\begin{equation}\label{Cauchy}
 \partial_t u = \Delta\,u +\beta\,u +\Psi_1(x)\,u^{-1}-\Psi_2(x)\,u^{-3},\quad u(x,0)=u_0(x)>0\,.
\end{equation}
Let $\mathcal{C}_{t}=F\times[0,t)$ be cylinder with the base $F$.
By~\cite[Theorem~4.51]{aub},
(\ref{Cauchy}) has a unique smooth solution in $\mathcal{C}_{t_0}$ for some $t_0>0$.
Let $S_t$ be a map which relates to each initial value $u_0\in C$ the value
of this solution at $t\in[0,t_0)$.
Since rhs of (\ref{Cauchy}) does not depend explicitly on $t$, the family $\{S_t\}$ has
the semigroup property, and it is a semigroup (i.e., $t_0=\infty$) when (\ref{Cauchy})
has a global solution for any $u_0(x)\in C$.

Let $\mathcal{H}:=-\Delta -\beta$ be a Schr\"{o}\-din\-ger operator with domain in~$H^2$,
and $\sigma(\mathcal{H})$ the spectrum.
One can add a real constant to $\beta$ such that $\mathcal{H}$ becomes invertible in
$L_2$ (e.g. $\lambda_0>0$) and $\mathcal{H}^{-1}$ is bounded in $L_2$.

\vskip1mm\noindent
\textbf{Theorem} (Elliptic regularity, see \cite{aub}).
\textit{If $\,0\notin\sigma(\mathcal{H})$, then for any integer $k\ge0$ we have}
 $\mathcal{H}^{-1}:\; H^k\rightarrow H^{k+2}$.

\vskip1mm

By the Elliptic regularity Theorem with $k=0$, we
have $\mathcal{H}^{-1}: L_2\rightarrow H^2$, and the embedding of
$H^2$ into $L_2$ is continuous and compact, see \cite{aub}.
Hence, the operator $\mathcal{H}^{-1}:\,L_2\rightarrow L_2$ is compact.
Thus, the spectrum $\sigma(\mathcal{H})$ is discrete, i.e.,
consists of an infinite sequence of real eigenvalues
$\lambda_0\le\lambda_1\le\ldots\le\lambda_j\le\ldots$ with finite multiplicities,
bounded from below and $\lim\nolimits_{j\to\infty}\lambda_j=\infty$.
One may fix in $L_2$ an orthonormal basis of eigenfunctions $\{e_j\}$, i.e., $\mathcal{H}(e_j)=\lambda_j e_j$.
Since the eigenvalue $\lambda_0$ is simple, its eigenfunction $e_0(x)$ (called the \textit{ground state})
can be chosen positive, see \cite[Prop.~3]{rz2012}.

\smallskip
 The following examples show us that (\ref{Cauchy}) may have

(i)~stationary (i.e., $t$-independent) solutions on a closed manifold $F$;

(ii)~attractors
(i.e., asymptotically stable stationary solutions) when $\beta<0$.

\begin{example}\label{Ex-2const}\rm
 Let $\beta$ and $\Psi_1>0$, $\Psi_2\ge0$ be real constants.
 Then (\ref{Cauchy}) is the Cauchy's problem for ODE
\begin{equation}\label{zermod}
 y\,'=f(y),\quad y(0)=y_0>0,\quad
 f(u):=\beta u +\Psi_1\,u^{-1} -\Psi_2\,u^{-3}.
\end{equation}

(a) Let $\beta<0$ and $\Psi_2>0$. Positive stationary (i.e., constant) solutions of (\ref{zermod})
are the roots of a biquadra\-tic equation $y^3f(y)=0$.
If $4\,|\beta|\Psi_2<\Psi_1^2$, then we have two positive solutions
$y_{1,2}=\big(\frac{\Psi_1\pm(\Psi_1^2-4\,|\beta|\,\Psi_2)^{1/2}}{2 |\beta|}\big)^{1/2}$
and $y_1>y_2$.
The linearization of (\ref{zermod}) at the point $y_k\ (k=1,2)$ is
 $v\,'=f^\prime(y_k)v$,
where
 $f^\prime(y_k)=-|\beta|\,(y^{-3}(y^2-y_1^2)(y^2-y_2^2))^\prime\vert_{\,y=y_k}$.
Hence, $f^\prime(y_1)<0$ and $f^\prime(y_2)>0$, and
$y_1$ is asymptotically stable, but $y_2$ is unstable.
If $4\,|\beta|\Psi_2=\Psi_1^2$, then (\ref{zermod}) has one
positive stationary solution, see Fig.~\ref{F-01}(a),
and has no stationary soluti\-ons if $4\,|\beta|\Psi_2>\Psi_1^2$.

\begin{figure}
\begin{tabular*}{\textwidth}{@{}l@{}l@{}r@{}}
\hskip5mm
\subfigure[\small $\beta\,y^4 +\Psi_1 y^{2} -\Psi_2$ with $\beta<0$ and $4\,|\beta|\,\Psi_2<\Psi_1^2$:
$y_1$~stable, $y_2$ unstable]{\includegraphics[width=0.35\textwidth]{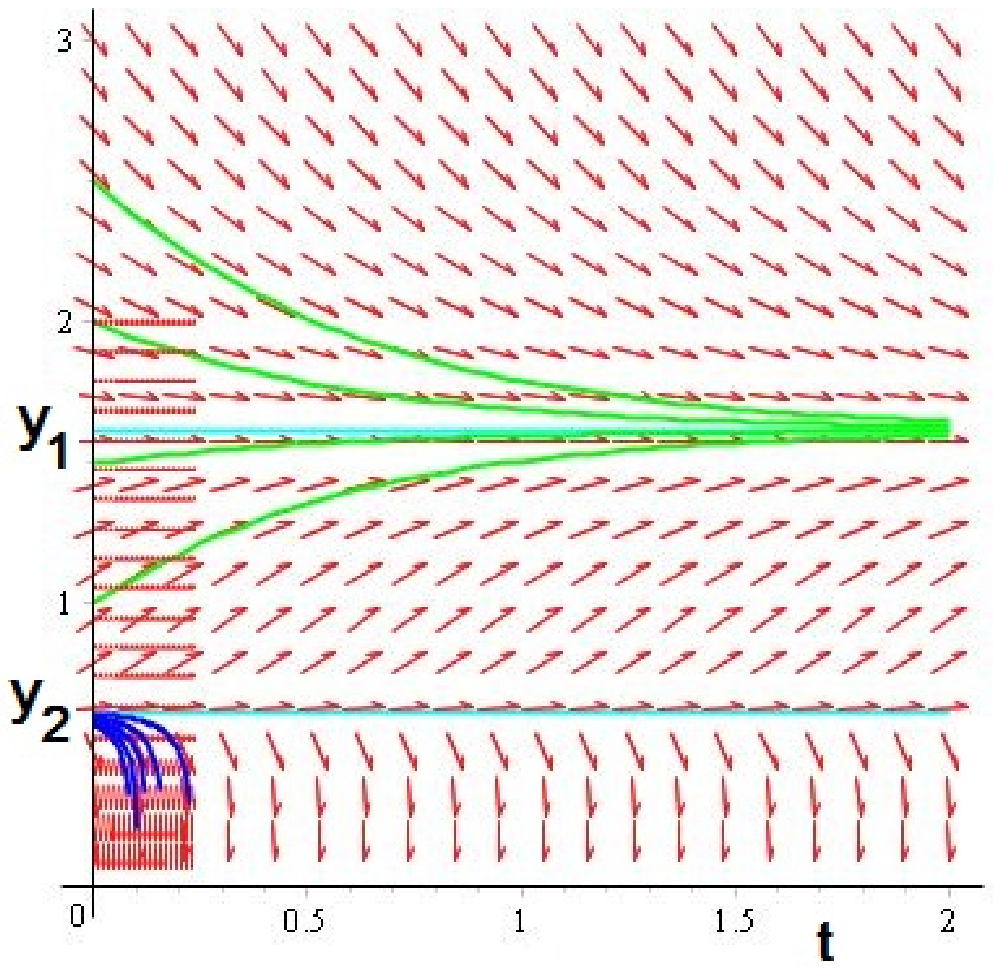}
}
\hskip10mm
\subfigure[\small
$\Psi_1>0$, $\Psi_2=0$   and $\beta<0$]{\includegraphics[width=0.43\textwidth]{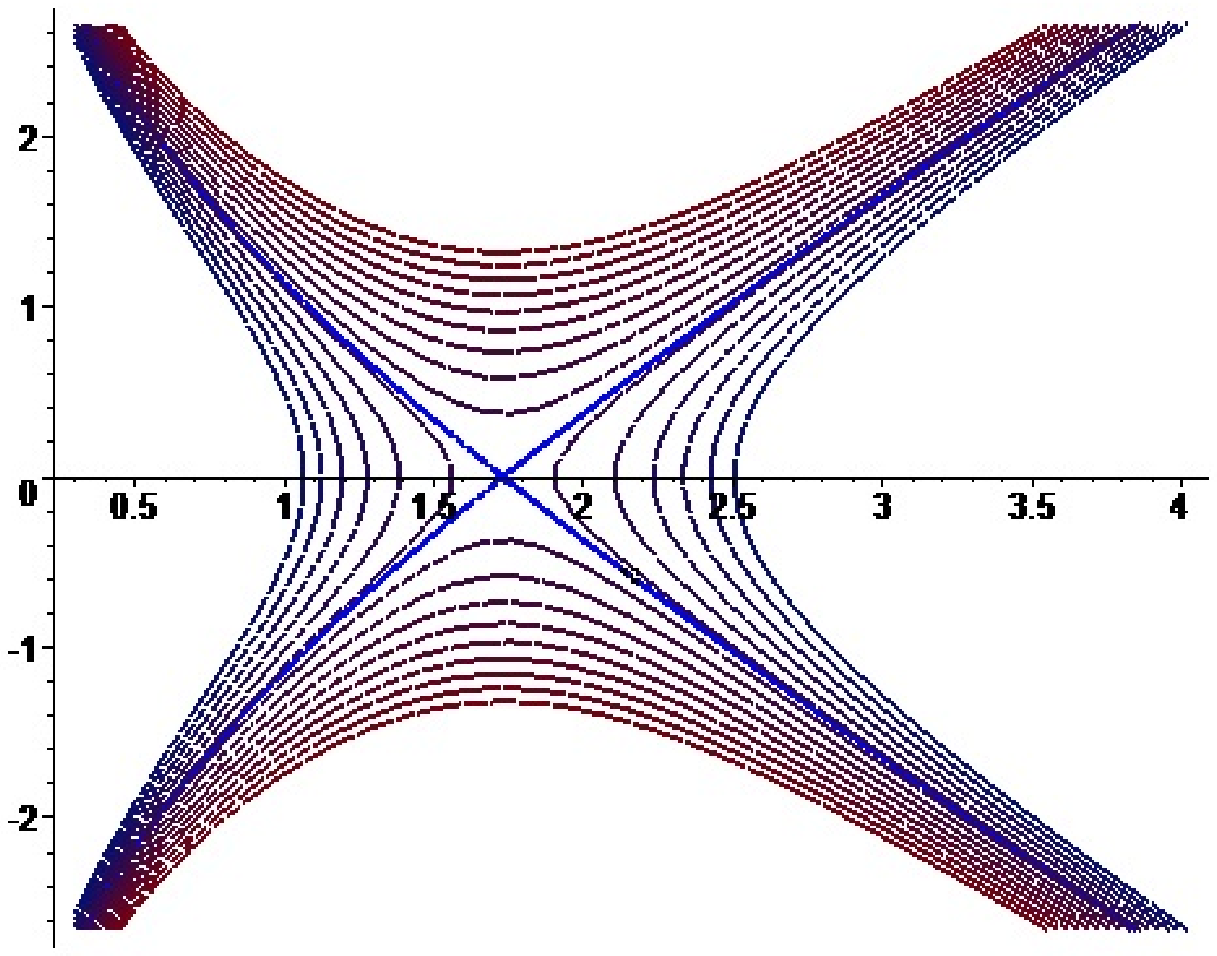}
}
\end{tabular*}
\caption{Example~\ref{Ex-2const}: the nonlinear heat equation}
\label{F-01}
\end{figure}

(b) Let $\beta>0$ and $\Psi_2>0$. Then the biquadratic equation $y^3f(y)=0$
has one positive root $y_1 = \big(\frac{-\Psi_1 +(\Psi_1^2+4\,\beta\,\Psi_2)^{1/2}}{2\,\beta}\big)^{1/2}$.
We find
\[
 f^\prime(y_1)=\beta\,
 \big(y^{-3}\big(y^2-y_1^2\big)
 \big(y^2 +{\Psi_2}/({\beta\,y_1^2})\big)\big)^\prime_{\,|\,y=y_1}>0;
\]
hence, $y_1$ is unstable. One may also show that in the case
$\beta=0$, (\ref{zermod}) has a unique positive stationary solution,
which is unstable.

(c) Let $\Psi_2=0$ and $\Psi_1>0$. Then $f(y)=\beta\,y+\Psi_1 y^{-1}$.
If $\beta\ge 0$, then there are no positive stationary solutions.
If $\beta<0$, then $f$ has one positive root $y_1=({\Psi_1/|\beta|})^{1/2}$.
Since
 $f^\prime(y_1)=-|\beta|\,(y^{-1}(y-y_1)(y+y_1))^\prime_{|\,y=y_1}<0$,
the solution $y_1$ is an attractor.
\end{example}

\begin{example}\label{Ex-2Bconst}\rm
 Let $F$ be a circle $S^1$ of length $l$.
Then (\ref{Cauchy}) is the Cauchy's problem
\begin{equation}\label{heat}
 u_{,t}=u_{,xx} +f(u),\quad u(x,0)= u_0(x)>0\quad (x\in S^1,\;t\ge 0).
\end{equation}
The stationary equation with $u(x)$ for (\ref{heat}) has the form
\begin{equation}\label{stateq2}
 u\,''+f(u)=0,\quad
 u(0)=u(l),\quad u\,'(0)=u\,'(l),\quad l>0.
\end{equation}
 Rewrite (\ref{stateq2}) as the dynamical system
\begin{equation}\label{Hamiltsyst}
 u\,'=v,\quad v\,'=-f(u)\quad (u>0).
\end{equation}
Periodic solutions of (\ref{stateq2}) correspond  to solutions of (\ref{Hamiltsyst}) with the same period.
The system (\ref{Hamiltsyst}) is Hamiltonian, since $\partial_u v=\partial_v f(u)$,
its {Hamiltonian} $\mathrm{H}(u,v)$ (the first integral) solves
 $\partial_u \mathrm{H}(u,v)=f(u),\ \partial_v \mathrm{H}(u,v)=v$.
 Then
 $\mathrm{H}(u,v)=\frac12(v^2+\beta u^2)+\Psi_1\ln u +\frac12\Psi_2\,u^{-2}$.
The trajectories of (\ref{Hamiltsyst}) belong to level lines of $\mathrm{H}(u,v)$.
Consider the~cases.

(a) Assume $\beta<0$. Then (\ref{Hamiltsyst}) has two fixed points:
$(y_i,0)\ (i=1,2)$~with $y_1>y_2$.
To clear up the type of fixed points, we linearize (\ref{Hamiltsyst}) at $(y_i,0)$,
\[
 \vec{\eta}\ '=A_i\,\vec\eta,\quad  A_i=\left(\begin{array}{ll} \quad 0& 1\\
 -f^\prime(y_i)& 0
\end{array}\right).
\]
Since $f^\prime(y_1)<0$ and $f^\prime(y_2)>0$,
the point $(y_1,\,0)$ is a ``saddle'' and $(y_2,\,0)$ is a ``center".
The separatrix is $\mathrm{H}(u,v)=\mathrm{H}(y_1,0)$, i.e., see Fig.~\ref{F-02}(a),
\[
 v^2=|\beta|(u^2-y_1^{\,2})-2\,\Psi_1\ln({u}/{y_1})-\Psi_2(u^{-2}-y_1^{\,-2}).
\]
The separatrix divides the half-plane $u>0$ into three simply connected areas.
Then $(y_2,\,0)$ is a~unique minimum point of $\mathrm{H}$ in
$D=\{(u,\,v):\ \mathrm{H}(u,v)<\mathrm{H}(y_1,0),\;0<u<y_1\}$.
The~phase portrait of (\ref{Hamiltsyst}) in $D$ consists of the cycles surrounding the fixed
point $(y_2,\,0)$, all correspond to non-constant solutions of (\ref{stateq2}) with various~$l$.
Other two areas do not contain cycles, since they have no fixed points.

(b) Assume $\beta\ge 0$. Then (\ref{Hamiltsyst}) has one fixed point
$(y_1,0)$ and $f^{\,\prime}(y_1)>0$. Hence, $(y_1,0)$ is a ``center".
Since $(y_1,\,0)$ is a unique minimum of $\mathrm{H}(u,v)$ in the semiplane $u>0$,
the phase portrait of (\ref{Hamiltsyst}) consists of
the cycles surrounding the fixed point $(y_1,\,0)$, all correspond to non-constant solutions of (\ref{stateq2}) with various $l$,
see Fig.~\ref{F-02}(b).

 For $\Psi_2=0$ and $\Psi_1>0$, the Hamiltonian of (\ref{Hamiltsyst})
is $\mathrm{H}(u,v)=\frac12 (v^2+\beta\,u^2)+\Psi_1\ln u$. Solving $\mathrm{H}(u,v)=C$ with respect to $v$ and
substituting to (\ref{Hamiltsyst})$_1$, we get
 $u^\prime=\sqrt{-\beta\,u^2-2\,\Psi_1\ln u+2\,C}$.
If $\beta\ge 0$, then (\ref{Hamiltsyst}) has no cycles (since it has no fixed points); hence, (\ref{stateq2})
has no solutions. If $\beta<0$, then the separatrix $\mathrm{H}(u,v)=\mathrm{H}(u_*,0)$, see Example~\ref{Ex-2const}(c),
is $v^2=|\beta|(u^2-u_*^2)-2\,\Psi_1\ln({u}/{u_*})$,
(\ref{Hamiltsyst}) has a unique fixed point $(u_*,0)$ which is a ``saddle".
The separatrix divides the half-plane $u>0$ into four simply connected areas with these lines,
see Fig.~\ref{F-01}(b).
Since each of these areas has no fixed points of (\ref{Hamiltsyst}),
the system has no cycles. Hence, $u_*$ is a unique solution of (\ref{stateq2}).

\begin{figure}
\begin{tabular*}{\textwidth}{@{}l@{}l@{}r@{}}
\hskip5mm
\subfigure[\small $\beta<0$]{\includegraphics[width=0.40\textwidth]{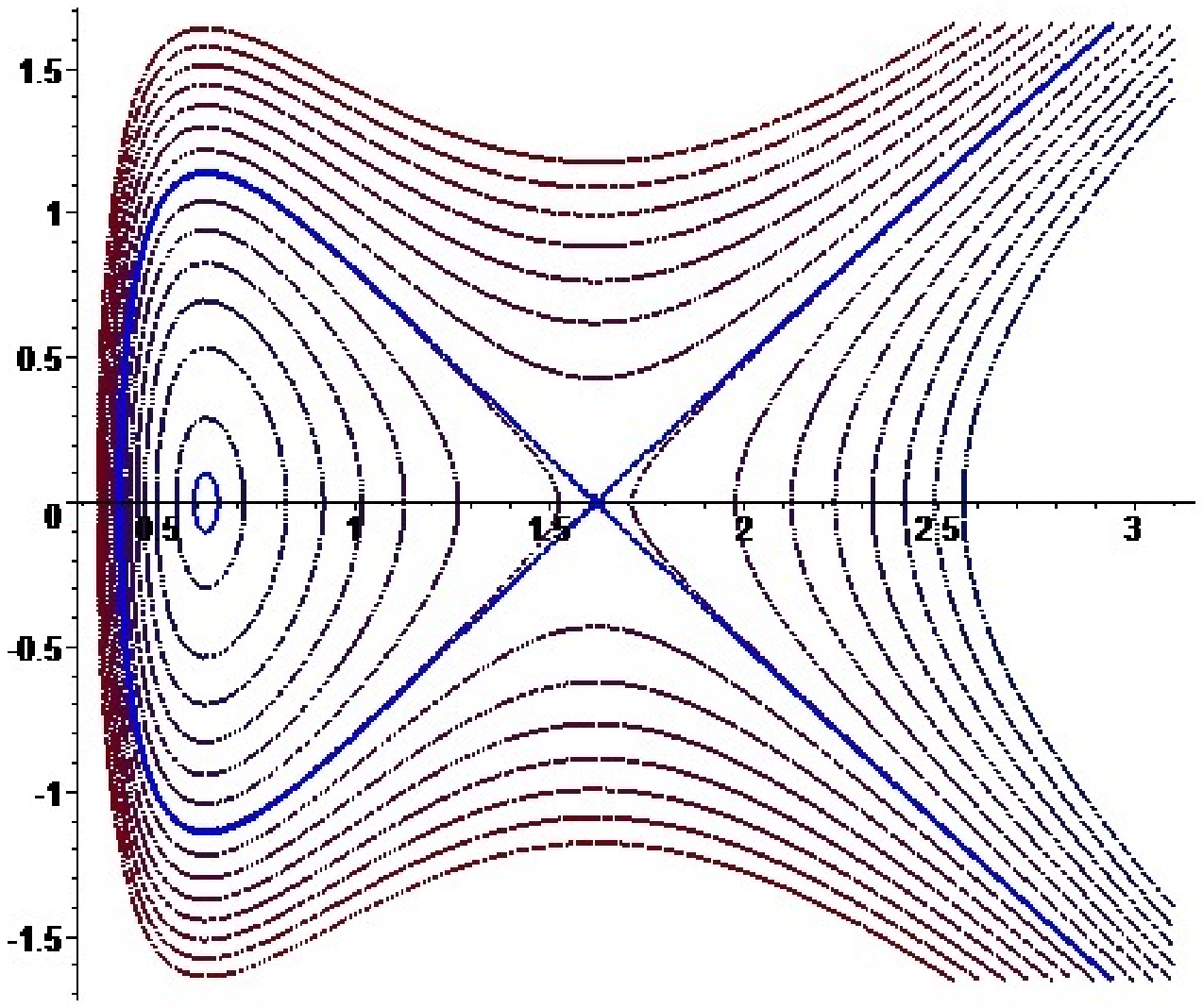}}
\hskip5mm
\subfigure[\small $\beta>0$]{\includegraphics[width=0.40\textwidth]{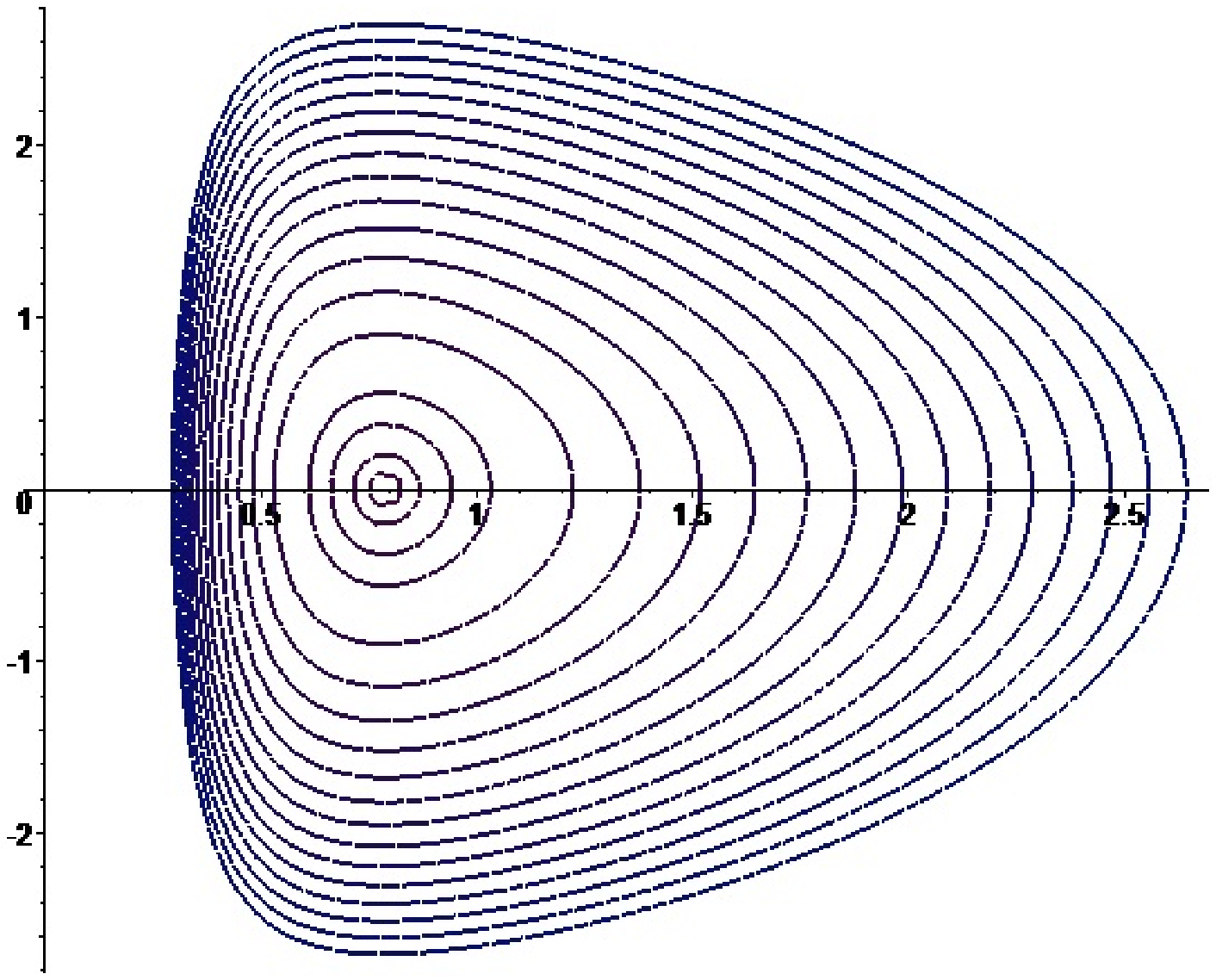}}
\end{tabular*}
\caption{Example~\ref{Ex-2Bconst}.}
\label{F-02}
\end{figure}

(c) Consider (\ref{stateq2}) for $\Psi_1=0$, $\Psi_2>0$ and $l=2\pi$.
Set $p=u\,'$ and represent $p=p(u)$ as a function of $u$. Then $u\,''={dp}/{du}$ and
\[
 (p^2)'=-2\,\beta\,u+2\,\Psi_2\,u^{-3}\ \Rightarrow\ (u\,')^2=C_1-\beta\,u^2-\Psi_2\,u^{-2}.
\]
After separation of variables and integration, we obtain
\[
  \ u=\!\left\{
  \begin{array}{cc}
    \sqrt{\frac{C_1}{2\beta}{+}\frac1{2\beta}\sqrt{C_1^2-4\beta\,\Psi_2}\sin(2\sqrt{\beta}(x{+}C_2))},
    \ \ (C_1^2\ge 4\beta\Psi_2) & {\rm for}\ \beta>0, \\
    \sqrt{-\frac{C_1}{2|\beta|}+\frac1{2|\beta|}\sqrt{C_1^2+4|\beta|\Psi_2}\cosh(2\sqrt{|\beta|}(x+C_2))}
    & {\rm for}\ \beta<0, \\
    \sqrt{{\Psi_2}/{C_1}+C_1(x+C_2)^2} & {\rm for}\ \beta=0.
  \end{array}\right.
\]
Hence, for $\beta\le0$, (\ref{stateq2}) has no positive solutions,
while for $\beta>0$ the solution is $2\pi$-periodic and positive only if

\noindent\
 $\bullet$ $\beta\ne{n^2}/{4}\ (n\in\NN)$ and $C_1=2\,(\beta\,\Psi_2)^{1/2}$;
 a solution $u_*=(\Psi_2/\beta)^{1/4}$ is~unique, or

\noindent\
 $\bullet$ $\beta={n^2}/{4}\ (n\in\NN)$; the set of solutions forms a two-dimensional manifold
\[
 n\,u_0(C_1,C_2)=\big({2\,C_1}+2(C_1^2-n^2\Psi_2)^{1/2}\sin(n(x+C_2))\big)^{1/2}\,.
\]
\end{example}

\begin{figure}
\begin{tabular*}{\textwidth}{@{}l@{}l@{}r@{}}
\hskip15mm\subfigure
[\small \ $y_1^-\approx 3, y_2^-\approx 1, y_3^-\approx 1.6$ and $y_4^-\approx 2.4$]
{\includegraphics[width=0.42\textwidth]{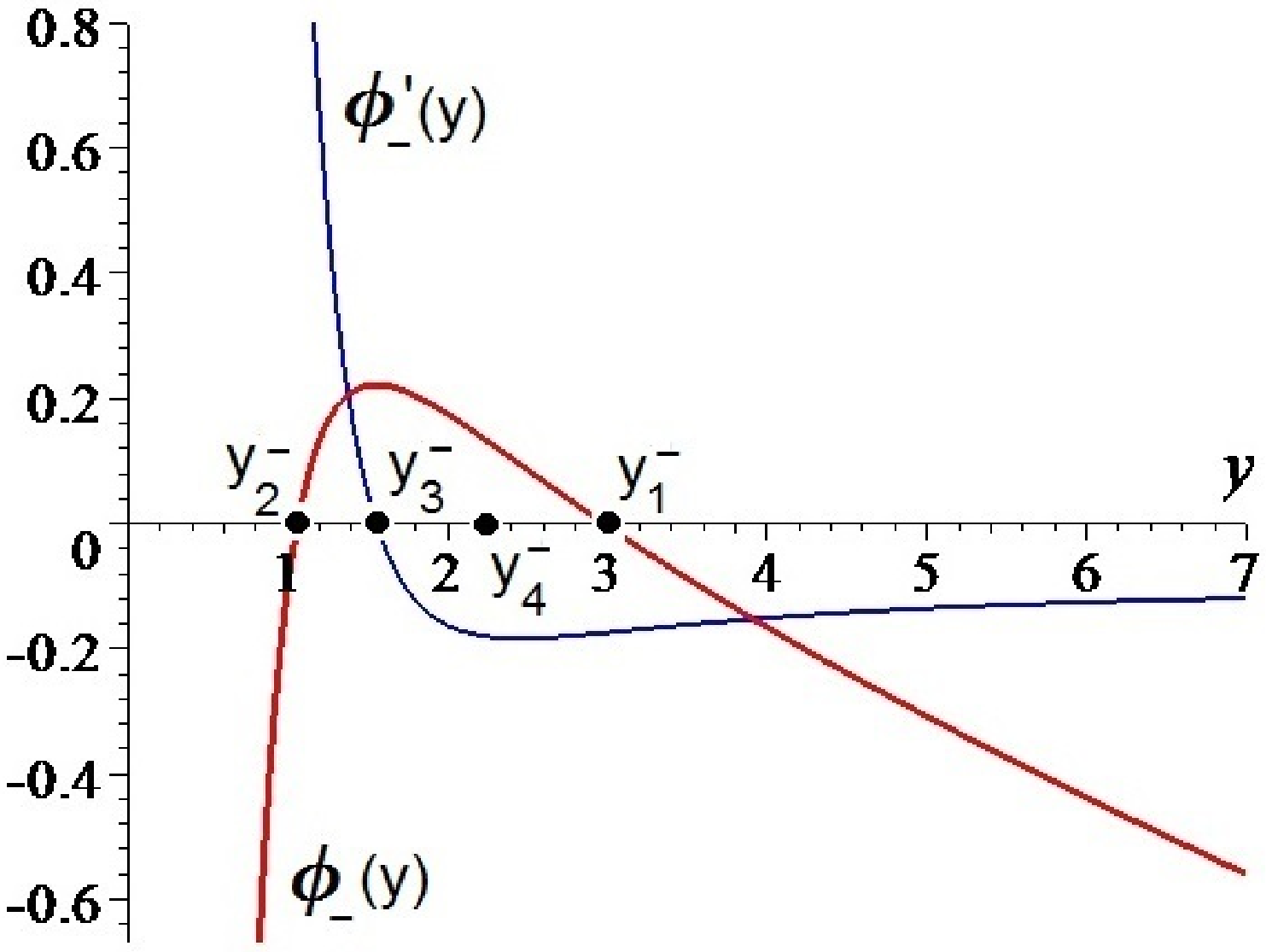}
}
\hskip5mm\subfigure
[\small $\mu^-(\sigma)$ for $0\le \sigma<y_1^- -y_3^-\approx 1.4$]
{\includegraphics[width=0.36\textwidth]{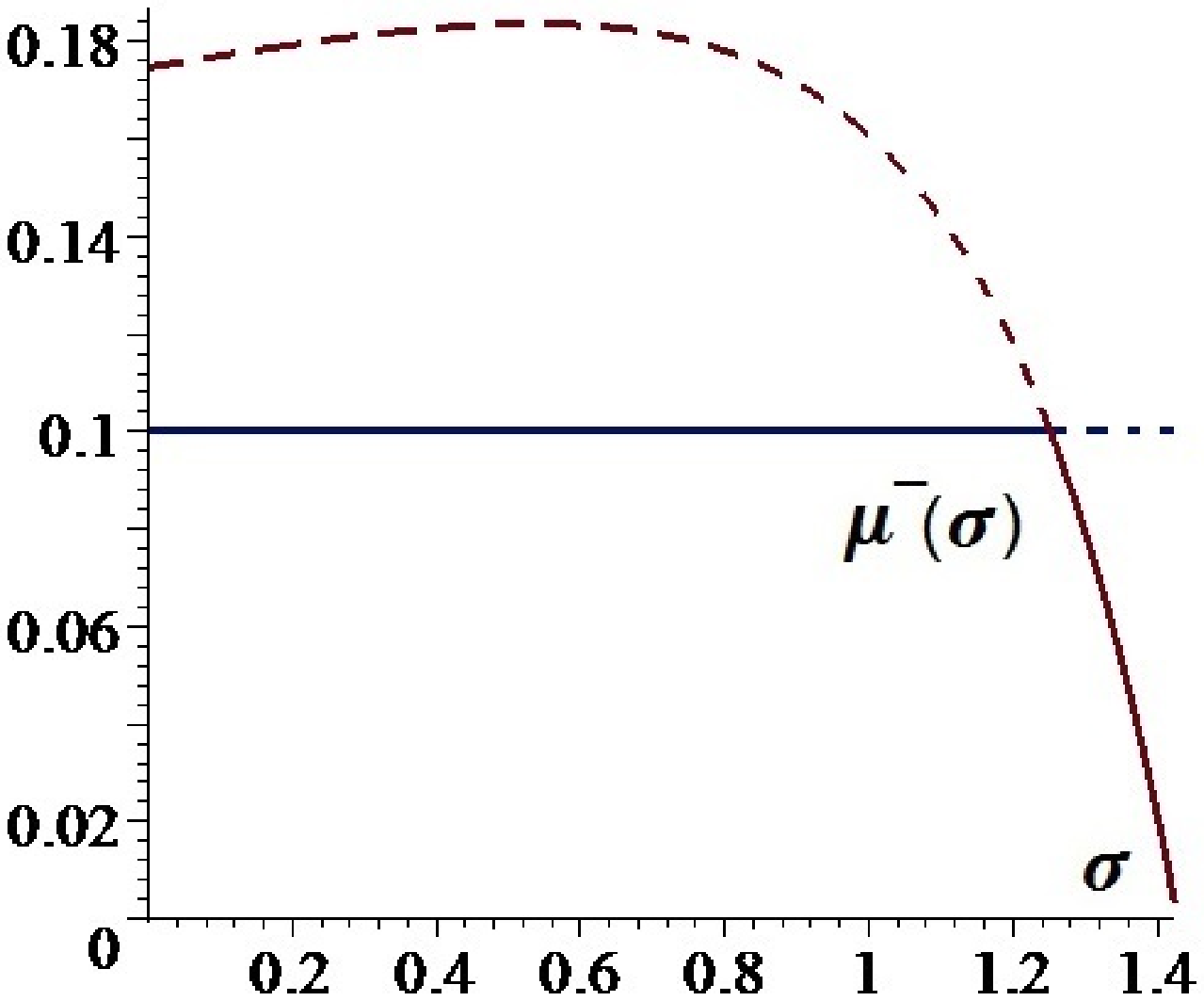}
}
\end{tabular*}
\vskip-2mm
\caption{Graphs of $\phi_-$, $\phi_-^\prime$ and $\mu^-$ for $\Psi_1=\Psi_2=1$ and $\lambda_0=0.1$.}
\label{fig:1}
\end{figure}

\textbf{2.2. Attractor of the nonlinear heat equation}.
Denote by  $\Psi_i^+=\max\nolimits_{\,F}\,(\Psi_i\,e_0^{-2i})$ and
 $\Psi_i^-=\min\nolimits_{\,F}\,(\Psi_i\,e_0^{-2i})$ for $i=1,2$.
Let $\Psi_2^+>0$ (the case of $\Psi_2^+=0$ is similar) and
\begin{equation}\label{E-lambda0-cond1}
 0<\lambda_0<{(\Psi_1^-)^2}/({4\,\Psi_2^+}).
\end{equation}
Each of the two functions of variable $y>0$,
\begin{equation}\label{E-phi-both}
 \phi_+(y)=-\lambda_0y +\Psi_1^+ y^{-1}\!-\Psi_2^- y^{-3},\quad
 \phi_-(y)=-\lambda_0y +\Psi_1^- y^{-1}\!-\Psi_2^+ y^{-3},
\end{equation}
has four real roots, two of which, $y_2^+<y_1^+$ and $y_2^-<y_1^-$, are positive.
Since $\phi_-(y)\le\phi_+(y)$ for $y>0$, we also have $y_1^- \le y_1^+$.
 Denote by
\begin{equation}\label{dfy3min}
 y_3^-=\Big(\frac{-\Psi_1^- +((\Psi_1^-)^2+12\,\Psi_2^+\lambda_0)^{1/2}}{2\,\lambda_0}\Big)^{1/2}
\end{equation}
a uni\-que positive root of $\phi_-^{\,\prime}(y)$. Cle\-arly, $y_3^-\in (y_2^-,\,y_1^-)$.
Notice that $\phi_-(y)>0$ for $y\in(y_2^-,\,y_1^-)$ and
$\phi_-(y)<0$ for $y\in(0,\,\infty)\setminus[y_2^-,\,y_1^-]$;
moreover, $\phi_-(y)$ increases in $(0,y_3^-)$ and decreases in $(y_3^-,\,\infty)$.
 The line $z=-\lambda_0\,y$ is asymptotic for the graph of $\phi_-(y)$ when $y\to\infty$,
 and $\lim\limits_{\,y\downarrow 0}\phi_-(y)=-\infty$.
 Next, $\phi_-^{\,\prime}(y)$ decreases in $(0,y_4^-)$ and increases in
$(y_4^-,\,\infty)$, where $y_4^-:=(6\,\Psi_2^+/\Psi_1^-)^{1/2}>y_3^-$,
and $\lim\limits_{\,y\rightarrow\infty}\phi_-^{\,\prime}(y)=-\lambda_0$,
see Fig.~\ref{fig:1}.
Hence,
\begin{equation}\label{E-mu-pm}
  \mu^-(\sigma):=-\sup\nolimits_{\,y\ge y_1^- -\sigma}\,\phi_-^\prime(y)
 =\min\{|\phi_-^\prime(y_1^- -\sigma)|,\,\lambda_0\}>0
\end{equation}
for $\sigma\in(0,\,y_1^- -y_3^-)$.
Similar properties have $y_3^+, y_4^+$ and $\mu^+(\sigma)$  defined for $\phi_+^\prime(y)$.

\begin{lemma}\label{lmstabODE}
Let $y(t)$ be a solution of Cauchy's problem
\begin{equation}\label{CauchODE}
 y\,'=\phi_-(y),\quad y(0)=y_0^->0.
\end{equation}

$(i)$ If $y_0^- > y_2^-$ then $\lim\nolimits_{\,t\rightarrow\infty}y(t)=y_1^-$.
Furthermore, if $y_0^-\in(y_2^-,y_1^-)$ then $y(t)$

\quad is increasing and if $y_0^->y_1^-$ then $y(t)$ is decreasing.

\smallskip
$(ii)$~If $y_0^-\ge y_1^- -\eps$ for some $\eps\in(0,\,y_1^- -y_3^-)$ then
\begin{equation}\label{eststab}
 |y(t)-y_1^-|\le|y_0^- - y_1^-|\,e^{-\mu^-(\eps)\,t}.
\end{equation}
Similar claims are valid for Cauchy's problem
$y\,'=\phi_+(y),\ y(0)=y_0^+>0$.
\end{lemma}

\noindent\textbf{Proof}. (i) Assume that $y_0^-\in(y_2^-,\,y_1^-)$.
Since $\phi_-(y)$ is positive in $(y_2^-,\,y_1^-)$, $y(t)$ is increasing.
The graph of $y(t)$ cannot intersect the graph of the stationary solution $y_1^-$;
hence, the solution $y(t)$ exists and is continuous on the whole $[0,\infty)$, and it is bounded there.
There exists $\lim\nolimits_{\,t\rightarrow\infty}y(t)$, which coincides with $y_1^-$,
since $y_1^-$ is a unique solution of $\phi(y)=0$ in $(y_2^-,\infty)$.
The~case $y_0^- > y_1^-$ is treated similarly.
Notice that if $y_0^- \in(y_2^-,\, y_1^-)$ then $y(t)$ is increasing,
and if $y_0^- > y_1^-$ then $y(t)$ is decreasing.

(ii) For $y_0^-\ge y_1^- -\eps$, where $\eps\in(0,\, y_1^- - y_3^-)$, denote $z(t)=y_1^- - y(t)$.
We obtain from (\ref{CauchODE}), using definition of $\mu^-(\eps)$ and the fact that $\phi_-(y_1^-)=0$,
\begin{eqnarray*}
 (z^2)'=2\,z\,z' = 2 z^2\int_0^1 \phi\,'_-(y+\tau z)\dtau \le -2\,\mu^-(\eps)\,z^2.
\end{eqnarray*}
This differential inequality implies (\ref{eststab}).
The~case $y_0^- > y_1^-$ is treated similarly.
\qed

\smallskip

Under assumption (\ref{E-lambda0-cond1}), define nonempty sets
${\mathcal U}_{\;2}^{\;\eps,\,\eta}\subset {\mathcal U}_{\;1}^{\;\eps}$, closed in $C$,
with $\eps\in(0,\,y_1^--y_3^-)$ and $\eta>0$ by
\[
 {\mathcal U}_{1}^{\eps}=\{u_0\in C:\ u_0/e_0\ge y_1^- -\eps\},\quad
 {\mathcal U}_{2}^{\eps,\eta}=\{u_0\in C:\ y_1^--\eps\le u_0/e_0 \le y_1^+ +\eta\}.
\]
Then, $\mathcal{U}_{\,1}^{\,\eps}\subset\mathcal{U}_{\,1}$,
where ${\mathcal U}_{\;1}\!=\!\{u_0\in C: u_0/e_0 > y_3^-\}$ is open in $C$.

\begin{proposition}\label{T-new1}
Let (\ref{E-lambda0-cond1}) hold. Then Cauchy's problem (\ref{Cauchy})
with $u_0\in {\mathcal U}_{\;1}^{\;\eps}$ for some $\eps\in(0,\,y_1^--y_3^-)$, admits a~unique global solution.
Furthermore, the sets ${\mathcal U}_{\;1}^{\;\eps}$ and ${\mathcal U}_{\;2}^{\;\eps,\eta}\ (\eta>0)$
are invariant for the semigroup of operators corresponding to (\ref{Cauchy})$_1$.
\end{proposition}

\noindent\textbf{Proof}.
Let $u(\,\cdot\,,t)\ (t\ge0)$ solve (\ref{Cauchy}) with $u_0\in\mathcal{U}_{\,1}^{\,\eps}$
for some $\eps\in(0,\,y_1^--y_3^-)$.
Substituting $u=e_0 w$ and using $\Delta e_0 +\beta e_0=-\lambda_0e_0$, yields the Cauchy's problem
\begin{equation}\label{Cauchw1}
\partial_t w =\Delta w+\,\<2\,\nabla\log e_0,\,\nabla w\> +f(w,\,\cdot\,),\quad
 w(\cdot\,,0) ={u_0}/{e_0}\ge y_1^--\eps,
\end{equation}
 for $w(x,t)$, where
\begin{equation}\label{E-fw}
 f(w,\,\cdot\,)=-\lambda_0\,w+(\Psi_1 e_0^{-2})\,{w^{-1}}-(\Psi_2 e_0^{-4})\,w^{-3}.
\end{equation}
From (\ref{Cauchw1}) and (\ref{E-phi-both}) we obtain the differential inequalities
\begin{equation}\label{difineq}
 \phi_-(w) \le \partial_t w - \Delta w - \<2\nabla\log e_0,\nabla w\> \le \phi_+(w).
\end{equation}
By Proposition~\ref{P-weak-max}, applied to the left inequality of (\ref{difineq}), and Lemma~\ref{lmstabODE},
in the maximal domain $D_M$ of the existence of the solution $w(x,t)$ of (\ref{Cauchw1}),
we obtain the inequality
\[
 w(\,\cdot\,,t)\ge y_1^- -\eps\,e^{-\mu^-(\eps)\,t}\ge y_1^--\eps > 0,
\]
 which implies that $w(x,t)$ cannot ``blowdown" to zero.
Since $\phi_+(w)\le{\Psi^+_1}{w^{-1}}$, from the right inequality of (\ref{difineq}),
applying Proposition~\ref{P-weak-max}, we obtain in $D_M$
\[
 w(\,\cdot\,,t)\le w_+(t)=(((u_0^+)^2-\Psi_1^+/\lambda_0)\,e^{-2\lambda_0 t}+\Psi_1^+/\lambda_0)^{1/2},
\]
where $w_+(t)$ solves the Cauchy's problem for ODE
\[
 {d w_+}/{dt} + \lambda_0 w_+ = {\Psi^+_1}{w_+^{-1}},\quad
 w_+(0)=u_0^+ := \max\nolimits_{\,F}\,(u_0/e_0)\,.
\]
By the above, the solution $u(x,t)$ of (\ref{Cauchy}) exists for all $(x,t)\in{\mathcal C}_\infty$,
and the set ${\mathcal U}_{\;1}^{\;\eps}$ is invariant for the semigroup of operators
$\,\mathcal{S}_t:u_0\mapsto u(\cdot\,,t)\ (t\ge 0)$ in $\mathcal{C}_\infty=F\times[0,\infty)$,
corresponding to (\ref{Cauchy})$_1$. Assuming $u_0\in{\mathcal U}_{\;1}^{\;\eps,\,\eta}$ and applying again
Proposition~\ref{P-weak-max} and Lemma~\ref{lmstabODE} to the right inequality of (\ref{difineq}),
we~get
\[
 w(\,\cdot\,,t)\le y_1^++\eta e^{-\mu^+(\sigma)\,t},\quad
 \sigma\in(0,\,y_1^+-y_3^+).
\]
Thus, $u(\cdot,t)\in{\mathcal U}_{\;1}^{\;\eps,\eta}\ (t>0)$.
Hence, also the set ${\mathcal U}_{\;2}^{\;\eps,\eta}$ is invariant for all~$\mathcal{S}_t$.
\qed

\begin{theorem}\label{thattract1}
$(i)$ If (\ref{E-lambda0-cond1}) holds then (\ref{Cauchy-stat1}) admits in $\,\mathcal{U}_{\,1}$ a unique solution $u_*$ (on $F$),
which is smooth; moreover, $u_*=\lim\nolimits_{\,t\to\infty} u(\cdot\,,t)$, where $u$ solves~(\ref{Cauchy}) with $u_0\in\mathcal{U}_{\,1}$, and $y_1^-\le u_*/e_0\le y_1^+$.
Furthermore, for any $\eps\in(0,\,y_1^- -y_3^-)$, the set $\mathcal{U}_{\,1}^{\,\eps}$
is attracted by (\ref{Cauchy})$_1$ exponentially fast to the point $u_*$ in $C$-norm:
\begin{equation}\label{rateconv1}
 \|u(\cdot\,,t)-u_*\|_{C} \le \delta_{e_0}^{-1}\,e^{-\mu^-(\eps)\,t}\| u_0-u_*\|_{C}
 \quad(t>0,\ u_0\in\mathcal{U}_{\,1}^{\,\eps}).
\end{equation}
$(ii)$
If $\beta,\Psi_1,\Psi_2$ are smooth functions on the product $F\times\RR^n$
with a smooth leafwise metric $g(\cdot,q)$ and (\ref{E-lambda0-cond1}) holds for any leaf $F\times\{q\}\ (q\in\RR^n)$
then the leafwise solution $u_*$ of (\ref{Cauchy-stat1}) is smooth on $F\times\RR^n$.
\end{theorem}

\noindent\textbf{Proof}.
(i) By Proposition~\ref{T-new1}, the set $\mathcal{U}_{\,1}^{\,\eps}$ is invariant for
the semigroup of operators $\,\mathcal{S}_t:u_0\to u(\cdot\,,t)\ (t\ge 0)$ corresponding to (\ref{Cauchy})$_1$, i.e.,
${\mathcal S}_t\big(\mathcal{U}_{\,1}^{\,\eps}\big)\subseteq\mathcal{U}_{\,1}^{\,\eps}$ for $t\ge 0$.
Take $u_i^0\in\mathcal{U}_{\,1}^{\,\eps}\ (i=1,2)$ and denote by
\[
 u_i(\cdot\,,t)={\mathcal S}_t(u_i^0),\quad
 w_i(\,\cdot,t)=u_i(\,\cdot,t)/e_0,\quad
 w_i^0=u_i^0/e_0.
\]
  Using (\ref{Cauchw1}) and the equalities
\[
 2\,\bar w\,\Delta\,\bar w=\Delta(\bar w^2)-2\,\|\nabla\,\bar w\|^2,\quad
 \nabla\,(\bar w^2)=2\,\bar w\nabla\,\bar w
\]
with $\bar w=w_2-w_1$, we obtain
\begin{eqnarray*}
 &&\partial_t\big((w_2-w_1)^2\big)=2\,(w_2-w_1)\,\partial_t(w_2-w_1)\le\Delta\big((w_2-w_1)^2\big)\\
 &&+\,\<2\,\nabla\log e_0,\,\nabla (w_2-w_1)^2\> +2\,(f(w_2,\,\cdot\,)-f(w_1,\,\cdot\,))(w_2-w_1).
\end{eqnarray*}
We estimate the last term, using $w_i\ge y_1^- -\eps>y_3^-\ (i=1,2)$, (\ref{E-mu-pm}) and (\ref{E-fw}),
\begin{eqnarray*}
 &&(f(w_2,\,\cdot\,)-f(w_1,\,\cdot\,))(w_2-w_1)\\
 &&=(w_2-w_1)^2\!\int_0^1\partial_w f(w_1+\tau(w_2-w_1),\,\cdot\,)\dtau
   \le-\mu^-(\eps)(w_2-w_1)^2.
\end{eqnarray*}
Thus, the function $v=(w_2-w_1)^2$ satisfies the differential inequality
\[
 \partial_t v \le\Delta\,v +\<2\,\nabla\log e_0,\ \nabla v\> -2\,\mu^-(\eps)\,v.
\]
By Proposition~\ref{P-weak-max},
$v(\,\cdot\,,t)\le v_+(t)$, where $v_+(t)$ solves the Cauchy's problem for ODE:
\[
 v\,_+'=-2\,\mu^-(\eps)\,v_+(t),\quad v_+(0)=\| w_2^0-w_1^0\|_C^2.
\]
Thus,
\begin{eqnarray*}
\nonumber
 &&\|\mathcal{S}_t(u_2^0)-\mathcal{S}_t(u_1^0)\|_C
 \le \| w_2(\,\cdot\,,t)-w_1(\,\cdot\,,t)\|_C \cdot\max\nolimits_{\,F}e_0\\
 &&\le e^{-\mu^-(\eps)\,t}\,\| w_2^0-w_1^0\|_C\cdot\max\nolimits_{\,F}e_0
 \le \delta^{-1}_{e_0} e^{-\mu^-(\eps)\,t}\,\| u_2^0-u_1^0\|_C\,,
\end{eqnarray*}
i.e., the operators $\mathcal{S}_t\ (t\ge 0)$ corresponding to (\ref{Cauchy}) satisfy in
$\mathcal{U}_{\,1}^{\,\eps}$ the Lipschitz condition with respect to $C$-norm
with the Lipschitz constant $\delta^{\,-1}_{e_0} e^{\,-\mu^-(\eps)\,t}$.

By Proposition \ref{T-new1}, for any $t\ge 0$
the operator ${\mathcal S}_t$ for (\ref{Cauchy}) maps the set
$\mathcal{U}_{\,1}^{\,\eps}$, which is closed in $C$, into itself,
and for $t>\frac1{\mu^-(\eps)}\ln \delta^{-1}_{e_0}$ it is a contraction there.
Since all operators ${\mathcal S}_t$ commute one with another, they
have a unique common fixed point $u_*$ in $\mathcal{U}_{\,1}^{\,\eps}$.
Since $\eps\in(0,\,y_1^- -y_3^-)$ is arbitrary,
$u_*$ is a unique common fixed point of all ${\mathcal S}_t$ in the set $\mathcal{U}_{\,1}$.
For any $u_0\in\mathcal{U}_{\,1}^{\,\eps}$ and $t\ge 0$, (\ref{rateconv1}) holds.
Thus, $u_*\in C$ is a generalized solution of (\ref{Cauchy-stat1}).
By the Elliptic Regularity Theorem, $u_*\in C^\infty$ and it is a classical solution.
By Proposition~\ref{T-new1}, $\mathcal{U}^{\,\eps,\,\eta}_{\,2}\subset\mathcal{U}^{\,\eps}_{\,1}$
is also ${\mathcal S}_t$-invariant, hence $u_*\in \mathcal{U}_{\,2}^{\,\eps,\,\eta}$.
Since $\eps\in(0,\,y_1^- -y_3^-)$ and $\eta>0$ are arbitrary, we get $y_1^-\le u_*/e_0\le y_1^+$.

Notice that if the functions $\Psi_1$ and $\Psi_2$ are constant
then $\phi_+=\phi_-$, see (\ref{E-phi-both}); in this case, $u_*/e_0=y_1^+=y_1^-$ is constant, too.

\smallskip
(ii) Let $e_0(x,q)>0$ be the normalized eigenfunction for
the minimal eigenvalue $\lambda_0(q)$ of the operator $\mathcal{H}_q=-\Delta-\beta(x,q)$.
By Theorem~\ref{contdeponqeigfunct} (in Sect.~\ref{sec:app}),
$\lambda_0\in C^\infty(\RR^n)$ and $e_0\in C^\infty(F\times\RR^n)$,
hence $y_3^-$, defined by (\ref{dfy3min}), smoothly depends on $q$.
As we have proved in (i), for any $q\in\RR^n$ the stationary equation, see also (\ref{Cauchy-stat1}),
\begin{equation}\label{nonlinstatpar1}
 \Delta_q\,u + f(u,x,q)=0,
\end{equation}
with $f(u,x,q)=\beta(x,q)u+\Psi_1(x,q)\,u^{-1}-\Psi_2(x,q)\,u^{-3}$ has a unique solution $u_*(x,q)$
in the open set ${\mathcal U}_{\;1}(q)=\{u_0\in C(F\times\RR^n):\, u_0/e_0(\cdot,q)>y_3^-(q)\}$.

Since $y_3^-(q)$ and $e_0(x,q)$ are continuous,
there exist open neighborhoods $U_*\subseteq C^{k+2,\alpha}$ of $u_*(x,0)$ and $V_0\subset\RR^n$ of $0$ such that
\begin{equation}\label{inclneighb}
  U_*\subseteq{\mathcal U}_{\;1}(q)\quad \forall\, q\in V_0.
\end{equation}
We claim that all eigenvalues of the linear operator
 $\mathcal{H}_*=-\Delta_{0}-\partial_u\,f(u_*(x,0),x,0)$,
acting in $L_2$ with the domain $H^2$, are positive.
To show this, observe that $y_1^-(0)\le u_*(\cdot,0)/e_0(\cdot,0)\le y_1^+(0)$.
Let $\tilde u(x,t)$ be a solution of Cauchy's problem for the evolution equation
\begin{equation}\label{Cauchstar}
 \partial_t\tilde u=-\mathcal{H}_*(\tilde u),\quad \tilde u(x,0)=\tilde u_0(x)\in C.
\end{equation}
Using the same arguments as in the proof of (i), we obtain that the function
$v(x,t)={\tilde u^{\,2}(x,t)}\,{e^{-2}_0(x,0)}$ obeys the
differential inequality with $\mu^-_0=\min\{|\phi_-^\prime(y_1^-)|,\,\lambda_0\}>0$:
\[
 \partial_t v \le\Delta_{0}\,v +\<2\,\nabla\log e_0(\cdot,0),\ \nabla v\> -2\,\mu^-_0\,v\,.
\]
 By Proposition~\ref{P-weak-max}, $v(\,\cdot\,,t)\le v_+(t)$,
where $v_+(t)$ solves the Cauchy's problem for ODE
\[
 v\,_+'=-2\,\mu^-_0\,v_+,\quad v_+(0)=\|{\tilde u_0}/{e_0(\cdot,0)}\|_C^2\,;
\]
moreover, for any $\tilde u_0\in C$ the function $\tilde u(x,t)$
tends to $0$ exponentially fast, as $t\rightarrow\infty$.
On the other hand, if $\tilde\lambda_\nu$ is any eigenvalue of $\mathcal{H}_*$
and $\tilde e_\nu(x)>0$ the corresponding normalized eigenfunction
then $\tilde u=e^{-\tilde\lambda_\nu t}\tilde e_\nu$ solves (\ref{Cauchstar})
with $\tilde u_0(x)=\tilde e_\nu(x)$.
Thus, $\tilde\lambda_\nu>0$ that completes the proof of the claim.

Using Theorem~\ref{prsmoothparstat}, we conclude that
for any integers $k\ge 0$ and $l\ge 1$ we can restrict the neighborhoods $U_*$ of $u_*(x,0)$ and $V_0$ of $0$
in such a way that for any $q\in V_0$ there exists in $U_*$ a unique solution
$\tilde u(x,q)$ of (\ref{nonlinstatpar1}) and the mapping
$q\rightarrow \tilde u(\cdot,q)$ belongs to class $C^l(V_0,U_*)$.
In view of (\ref{inclneighb}), $\tilde u(\cdot,q)=u_*(\cdot,q)$ for any $q\in V_0$.
\qed

\section{Appendix: Elliptic equation with parameter}
\label{sec:app}

Let $F\times\RR^n$ be the product with a compact leaf $F$, and $g(\cdot,q)$ a leafwise Riemannian
metric (i.e., on $F_q=F\times\{q\}$ for $q\in\RR^n$) such that
the volume form of leaves $d\,{\rm vol}_F=|g|^{1/2}\,{\rm d} x$ depends on $x\in F$ only.
 The Laplacian in a local chart $(U,x)$ on $(F,g_{|\,\calf})$ is written as
 $\Delta\,u = \nabla_i(g^{ij}\,\nabla_j\,u) = |g|^{-1/2}\partial_{i}(|g|^{1/2} g^{ij} \partial_j\,u)$,
 see \cite{aub} with opposite sign.

This defines an elliptic operator $-\Delta_q$, where $q\in\RR^n$ is a parameter and $\Delta_0=\Delta$,
\begin{equation}\label{ellop}
 \Delta_q = g^{ij}(x,q) \partial^2_{ij} + b^j(x,q)\partial_j\,.
\end{equation}
Here $b^j=|g|^{-1/2}\partial_{i}(|g|^{1/2} g^{ij})$ are smooth functions in $U\times\RR^n$.

The Schr\"odinger operator $\mathcal{H}_q=-\Delta_q-\beta(x,q)$ acts in the Hilbert space $L_2$ with the domain $H^2$.
 Denote ${\mathcal H}_k={\mathcal H}_{\,|\,H^{k+2}}$
and $\mathcal{H}_{q,k}=(\mathcal{H}_{q})_{\,|\,H^{k+2}}$ for any $q\in\RR^n$.
Consider also products $\mathbb{B}=L_2\times\RR^n$ and $\mathbb{B}_k=H^k\times\RR^n$
(trivial vector bundles over~$\RR^n$).

\smallskip
\textbf{3.1. The Schr\"odinger operator.}
If $B$ and $C$ are Banach spaces with norms $\|\cdot\|_B$ and $\|\cdot\|_C$, denote by $\mathfrak{B}^r(B,C)$ the Banach space
of all bounded $r$-linear operators $A:\,\prod_{i=1}^r B\rightarrow C$ with the norm
 $\|A\|_{\mathfrak{B}^r(B,C)}=\sup_{v_1,\dots, v_r\in B\setminus{0}}\frac{\|A(v_1,\dots,v_r)\|_C}{\|v_1\|_B\cdot\ldots\cdot\|v_r\|_B}$.
If~$r=1$, we shall write $\mathfrak{B}(B,C)$ and  $A(\cdot)$,
and if $B=C$ we shall write $\mathfrak{B}^r(B)$ and $\mathfrak{B}(B)$, respectively.
 If $M$ is a $k$-regular manifold or an open neighborhood of the origin in a real Banach space, and $N$ is a real Banach space,
we denote by $C^k(M,N)$ $(k\ge 1)$ the Banach space of all $C^k$-regular functions $f:\,M\rightarrow N$,
for which the following norm is finite:
\[
 \|f\|_{C^k(M,N)}=\sup\nolimits_{\,x\in M}\max\{\|f(x)\|_N,\,\max\nolimits_{1\le j\le k}\|d^j f(x)\|_{\mathfrak{B}^j(T_x M,\,N)}\}.
\]
 We~shall use the simplified version of the Banach's Closed Graph Theorem:

\textit{If a linear operator $A:\,B\rightarrow C$
(of Banach spaces $B$ and $C$) is bijective and boun\-ded, then its  inverse  $A^{-1}:\, C\rightarrow B$  is also bounded}.

\begin{lemma}\label{lmschrhilb}
Let $\beta\in C^\infty$ and $\mu<-\max_{x\,\in F}\,\beta(x)$. Then

$(i)$ ${\mathcal H}-\mu$ acts from $H^2$ into $L_2$, it is continuously inver\-tible and the inverse operator
$({\mathcal H}-\mu)^{-1}:\;L_2\rightarrow L_2$ is compact;

$(ii)$ for any $k\in\NN$ the operator ${\mathcal H}_k-\mu$ acts from $H^{k+2}$ into $H^k$,
it is continuously invertible and $({\mathcal H}_k-\mu)^{-1}:\;H^k\rightarrow H^k$ is compact;

$(iii)$ for any integer $k\ge0$ the spectrum of ${\mathcal H}_k$,
acting in $H^k$ with the domain $H^{k+2}$, is discrete, and it coincides with the spectrum
$\sigma({\mathcal H})$;

$(iv)$ for any integer $k\ge0$ and $\lambda\notin\sigma({\mathcal H})$ we have
\begin{eqnarray}\label{actRlambda}
 R_\lambda({\mathcal H}_k)=({\mathcal H}_k-\lambda)^{-1}\in\mathfrak{B}(H^k,\;H^{k+2}),\\
\label{mmbRlambda}
 (\lambda\rightarrow R_\lambda({\mathcal H}_k))\in C(\CC\setminus\sigma({\mathcal H}),\;\mathfrak{B}(H^k,\;H^{k+2})).
\end{eqnarray}
\end{lemma}

\noindent\textbf{Proof}.
(i) Clearly, there exists $C>0$ such that for any $u\in H^2$ we have
\begin{equation}\label{estnrmH}
 \Vert{\mathcal H}(u)-\mu u\Vert_{L_2}\le\Vert\Delta u-\mu\Vert_{L_2}
 +\Vert (\beta(x)+\mu) u\Vert_{L_2}\le C\Vert u\Vert_{H^2}.
\end{equation}
Thus, ${\mathcal H}-\mu\in\mathfrak{B}(H^2,\,L_2)$.
On the other hand, since $\beta(x)+\mu<0$, the operator ${\mathcal H}-\mu$ is positive definite as
acting in $L_2$ with the domain $H^2$;
hence, it is continuously invertible, i.e., $\mu\notin\sigma({\cal H})$.
Thus, ${\mathcal H}-\mu$ maps $H^2$ injectively onto~$L_2$. By~the Banach's Closed Graph Theorem,
${\mathcal H}-\mu:\,H^2\rightarrow L_2$ has the bounded inverse.
Compactness of the embedding $H^2\hookrightarrow L_2$ yields
compactness of $({\mathcal H}-\mu)^{-1}:\;H^k\rightarrow H^k$.

(ii) For any $u\in H^{k+2}$ a similar to (\ref{estnrmH}) estimate holds
with $H^k$ and $H^{k+2}$ instead of $L_2$ and $H^2$, respectively.
Hence,
 ${\mathcal H}_k-\mu\in\mathfrak{B}(H^{k+2},\,H^k)$.
By (i), ${\mathcal H}-\mu:\,H^2\rightarrow L_2$ is bijective,
and by Elliptic Regularity Theorem, $({\mathcal H}-\mu)^{-1}(H^k)\subseteq H^{k+2}$ holds.
Hence, ${\mathcal H}_k-\mu:\,H^{k+2}\rightarrow H^k$ is bijective, too.
Thus, the Banach's Closed Graph Theorem and the compactness of the embedding
$H^{k+2}\hookrightarrow H^k$ complete the proof of~(ii).

(iii) Since for any integer $k\ge0$ the operator
$({\mathcal H}_k-\mu)^{-1}:\;H^k\rightarrow H^k$ is compact,
its spectrum consists of $\nu=0$ and a countable number of
non-zero eigenvalues $\nu_n$ of finite multiplicity, which can accumulate only at $\nu=0$.
Hence, the spectrum of ${\mathcal H}_k$ on $H_k$ (with domain in $H^{k+2}$) is discrete.
Since ${\mathcal H}$ extends ${\mathcal H}_k$,
each eigenfunction of ${\mathcal H}_k$ is an eigenfunction of ${\mathcal H}$;
hence, $\sigma({\mathcal H}_k)\subseteq\sigma({\mathcal H})$.
 To show the opposite inclusion, let $e(x)$ be an eigenfunction
of ${\mathcal H}$, related to $\lambda\in\sigma({\mathcal H})$ (hence, $e\in H^2$).
The obvious equa\-lity $e=(\lambda-\nu)^j({\mathcal H}-\mu)^{-j}e$ is valid for any $j\in\NN$.
Applying (several times) the Elliptic Regularity Theorem, from $e\in H^{2}$ we obtain that $e\in H^{k+2}$.

(iv) Similarly to the proof of (ii), we obtain inclusion (\ref{actRlambda}).
Let us prove (\ref{mmbRlambda}). By~the well-known property of the resolvent, we have
\[
 (\lambda\rightarrow R_\lambda({\mathcal H}_k))\in C(\CC\setminus\sigma({\mathcal H}_k),\;\mathfrak{B}(H^k)).
\]
Take an arbitrary $\lambda_0\notin\sigma({\mathcal H})$ and choose $\delta>0$ such that the set
 $D_\delta=\{\lambda\in\CC:\;|\lambda-\lambda_0|\le\delta\}$
does not intersect $\sigma({\mathcal H})$.
Using the resolvent identity, see \cite{Akh-Gl},
\[
 R_\lambda({\mathcal H}_k)-R_{\lambda_0}({\mathcal H}_k)
 =(\lambda-\lambda_0)R_{\lambda_0}({\mathcal H}_k)R_\lambda({\mathcal H}_k),
\]
 we have the following estimate for $\lambda\in D_\delta$:
\begin{eqnarray*}
&&
 \Vert R_\lambda({\mathcal H}_k)-R_{\lambda_0}({\mathcal H}_k)\Vert_{\mathfrak{B}(H^k,\;H^{k+2})}
 \le|\lambda-\lambda_0|\times\\
&&
 \times\Vert R_{\lambda_0}({\mathcal H}_k)\Vert_{\mathfrak{B}(H^k,\;H^{k+2})}
 \max\nolimits_{\,\lambda\in D_\delta}\Vert R_\lambda({\mathcal H}_k)\Vert_{\mathfrak{B}(H^k)},
\end{eqnarray*}
which implies the desired inclusion (\ref{mmbRlambda}).
\qed

\smallskip
\textbf{3.2. The ground state.}
We will show smooth dependence on $q$ of the least eigenvalue $\lambda(q)$ of $\mathcal{H}_q$
and of the corresponding normalized eigenfunction $e(x,q)>0$.

\begin{lemma}\label{lmsmoothellipt1}
If $\beta\in C^\infty(F\times\RR^n)$ then for any $l\in\NN$ and integer $k\ge 0$ the mapping
$\mathfrak{D}:(u,q)\mapsto\mathcal{H}_q(u)$ is a $C^l$-morphism $($of trivial vector bundles$)$ from
$\mathbb{B}_{k+2}$ into $\mathbb{B}_k$.
\end{lemma}

\noindent\textbf{Proof}.
Fix a finite atlas $\{(U_a,x_a)\}_{1\le a\le A}$ on $F$, and let $\{\rho_a(x)\}_{1\le a\le A}$
be a subordinated partition of unity.
Taking $u\in H^{k+2}$, $\;q,s\in\RR^n$ and using (\ref{ellop}), we obtain
\begin{eqnarray*}
 &&1/t(\mathfrak{D}(u,q+ts)-\mathfrak{D}(u_a,q))-\mathfrak{D}_1(u,q)\,s\\
 &&=- \big(1/t\int_0^t\partial_q g^{ij}(x,q+\tau s)\,s\dtau
 -\partial_q g^{ij}(x,q)\,s\big)\partial^2_{ij}\,u\\
 &&- \big(1/t\int_0^t\partial_q b^i(x,q+\tau s)\,s\dtau
 -\partial_q b^i(x,q)\,s\big)\,\partial_i u\\
&&-\big(1/t\int_0^t\partial_q \beta(x,q+\tau s)\,s\dtau
 -\partial_q \beta(x,q)\,s\big)\,u
\end{eqnarray*}
in a local chart, where
\[
 \mathfrak{D}_1(\cdot,q)\,s = - \partial_q g^{ij}(x,q)\,s\partial^2_{ij}
 -\partial_q b^i(x,q)\,s\partial_i-\partial_q \beta(x,q)\,s.
\]
Hence,
\begin{eqnarray*}
 &&\Vert 1/t(\mathfrak{D}(u,q+ts)-\mathfrak{D}(u,q))-\mathfrak{D}_1(u,q)\,s\,\Vert_{H^k}^2
 =\sum\nolimits_{a=1}^A\int_{U_a}\rho_a(x) \\
 && \sum\nolimits_{|\vec m|\le k}\big|\partial_{\vec m}^{|\vec m|}\big(1/t(\mathfrak{D}(u,q+ts)-\mathfrak{D}(u,q))
 -\mathfrak{D}_1(u,q)\,s\big)\big|^2\,{\rm d}x\\
 &&\le C\max_{i,j\in\{1,2,\dots,p\}}\,\max_{a\in\{1,2,\dots,A\}}\,\max_{\tau\in[0,t]}
  \big(\Vert\partial_q g^{ij}(x,q+\tau s)-\partial_q g^{ij}(x,q)\Vert_{C^k(U_a,\,\mathfrak{B}(\RR^n))}^2\\
 &&+\Vert\partial_q b^i(x,q+\tau s)-\partial_q b^i(x,q)\Vert_{C^k(U_a,\,\mathfrak{B}(\RR^n))}^2\\
 &&+\Vert\partial_q \beta(x,q+\tau s)-\partial_q \beta(x,q)\Vert_{C^k(U_a,\,\mathfrak{B}(\RR^n))}^2\big)
 |s|^2\,\Vert u\Vert_{H^{k+2}}
\end{eqnarray*}
holds for some $C>0$ that does not depend on $u$.
We conclude that $\mathfrak{D}:\,H^{k+2}\times\RR^n\rightarrow H^k$
has the partial Gateaux differential $\partial_q\mathfrak{D}(u,q)$
at each point $(u,q)$, and it is equal to $\mathfrak{D}_1(u,q)\,s$.
Similarly, for any $(u_k,q_k)\in H^{k+2}\times\RR^n\;(k=1,2)$ we obtain
\begin{eqnarray*}
 &&\Vert \partial_q\mathfrak{D}(u_1,q_1)\,s-\partial_q\mathfrak{D}(u_2,q_2)\,s\,\Vert_{H^k}\\
 &&\le C\,|s|^2\max_{i,j\in\{1,2,\dots,p\}}\,\max_{a\in\{1,2,\dots,A\}}
 \big((\Vert\partial_q g^{ij}(x,q_1)-\partial_q g^{ij}(x,q_2)\Vert^2_{C^k(U_a,\,\mathfrak{B}(\RR^n))}\\
 && +\Vert\partial_q b^i(x,q_1)-\partial_q b^i(x,q_2)\Vert_{C^k(U_a,\,\mathfrak{B}(\RR^n))}+
 \Vert\partial_q\beta(x,q_1)-\partial_q\beta(x,q_2)\Vert_{C^k(U_a,\,\mathfrak{B}(\RR^n))})\\
 &&\times\Vert u_1\Vert_{H^{k+2}}+
 (\Vert\partial_q g^{ij}(x,q_2)\Vert_{\mathfrak{B}(\RR^n)}
 +\Vert\partial_q b^i(x,q_2)\Vert_{C^k(U_a,\,\mathfrak{B}(\RR^n))}\\
 &&+\Vert\partial_q\beta(x,q_2)\Vert_{C^k(U_a,\,\mathfrak{B}(\RR^n))})
 \cdot\Vert u_1-u_2\Vert_{H^{k+2}}\big),
\end{eqnarray*}
and conclude that the partial differential $\partial_q\mathfrak{D}$ is continuous:
\[
 \partial_q \mathfrak{D}\in C(H^{k+2}\times\RR^n,\;\mathfrak{B}(\RR^n,\;H^k)).
\]
One may prove by induction that for any $l\in\NN$ the mapping
$\mathfrak{D}(\cdot,\cdot)$ has at any point $(u,q)\in H^{k+2}\times\RR^n$ the partial differential
of $l$-th order $\partial_q^l \mathfrak{D}(u,q)$, and it has the form
\[
 \partial_q^l\mathfrak{D}(u,q) = -\partial_q^l g^{ij}(x,q)\partial^2_{ij}\,u
 -\partial_q^l b^i(x,q)\partial_i\,u-\partial_q^l\beta(x,q)\,u
\]
in a local chart, and
 $\partial_q^l\mathfrak{D}\in C(H^{k+2}\times\RR^n,\;\mathfrak{B}^l(\RR^n,\;H^k))$.
Since $\mathfrak{D}(u,q)$ and $\partial_q^l\mathfrak{D}(u,q)$ are linear by $u$,
this differential is continuous and $\mathfrak{D}(u,q)$ has continuous differentials by $q$ and $u$ of any order.
\qed

\begin{lemma}\label{lmcontres}
Let $K$ be a compact set of $\,\CC\setminus\sigma(\mathcal{H}_{0,k})$ for some integer $k\ge 0$.
If $\beta\in C^\infty(F\times\RR^n)$,
then there is an open neighborhood $W\subseteq\RR^n$ of the origin such~that
\begin{equation}\label{inclK}
 K\subset\CC\setminus\sigma(\mathcal{H}_{q,k})\quad \forall\,q\in W,
\end{equation}
and the following inclusion holds for any $l\in\NN$:
\begin{equation}\label{inclRlambda}
 (\lambda\rightarrow R_\lambda(\mathcal{H}_{q,k}))
 \in C(K,\;C^l(W,\;\mathfrak{B}(H^k, H^{k+2}))).
\end{equation}
\end{lemma}

\noindent\textbf{Proof}. The following obvious representation holds for $\lambda\in K$:
\begin{equation}\label{rprHqmuId}
 \mathcal{H}_{q,k}-\lambda=(\id+L(q,\lambda))(\mathcal{H}_{0,k}-\lambda)\quad (q\in{W}),
\end{equation}
where
\begin{equation}\label{dfLqlam}
 L(q,\lambda):=(\mathcal{H}_{q,k}-\mathcal{H}_{0,k})R_\lambda(\mathcal{H}_{0,k}).
\end{equation}
Using Lemma~\ref{lmsmoothellipt1}, we get that for any integer $l\ge 0$
\begin{equation}\label{mmbshpHq}
 (q\rightarrow \mathcal{H}_{q,k})\in C^l({W},\;\mathfrak{B}(H^{k+2},\;H^k)).
\end{equation}
Taking into account Lemma~\ref{lmschrhilb}(iv), we have for any $(q,\lambda)\in{W}\times K$
\begin{eqnarray*}
 \Vert L(q,\lambda)\Vert_{\mathfrak{B}(H^k)}\le\Vert\mathcal{H}_{q,k}-\mathcal{H}_{0,k}
 \Vert_{\mathfrak{B}(H^{k+2},\;H^k)}
 \max_{\mu\in K}\Vert R_\mu(\mathcal{H}_{0,k})\Vert_{\mathfrak{B}(H^k,\;H^{k+2})}.
\end{eqnarray*}
Hence, and in view of (\ref{mmbshpHq}) with $l=0$,
there exists an open neighborhood $W\subset\RR^n$ of the origin such that
\begin{equation}\label{boundLqmu}
 \sup\nolimits_{\,(q,\lambda)\in W\times K}\Vert L(q,\lambda)\Vert_{\mathfrak{B}(H^k)}
 \le 1/2\,.
\end{equation}
Thus, for any $(q,\lambda)\in W\times K$ the operator
$\id+L(q,\lambda)\in\mathfrak{B}(H^k)$ is continuously invertible
and its inverse is expressed by the Neumann series $
(\id+L(q,\lambda))^{-1}=\sum_{j=0}^\infty(-L(q,\lambda))^j$ converging in the $\mathfrak{B}(H^k)$-norm.
In view of (\ref{rprHqmuId}), we conclude that (\ref{inclK}) is valid and for any $(q,\lambda)\in W\times K$ we have
\begin{equation}\label{rprresHq}
 R_\lambda(\mathcal{H}_{q,k})=R_\lambda(\mathcal{H}_{0,k})(\id+L(q,\lambda))^{-1}.
\end{equation}
Lemma~\ref{lmschrhilb}(iv) and (\ref{dfLqlam})--(\ref{mmbshpHq}) imply
$L(\cdot,\lambda)\in C^l(W,\mathfrak{B}(H^k))$ for $\lambda\in K$ and $l\in\ZZ$,
and using the resolvent identity, we obtain
$L(q,\lambda)-L(q,\mu)=(\lambda-\mu)L(q,\mu)R_\lambda(\mathcal{H}_{0,k})$ for $\lambda,\mu\in K$.
Hence,
\[
 \Vert L(\cdot,\lambda)-L(\cdot,\mu)\Vert_{C^l(W,\mathfrak{B}(H^k))}
 \le |\lambda-\mu|{\cdot}\Vert L(\cdot,\mu)\Vert_{C^l(W,\mathfrak{B}(H^k))}
 \max_{\nu\in K}\Vert R_\nu(\mathcal{H}_{0,k})\Vert_{\mathfrak{B}(H^k)}.
\]
This estimate implies
\begin{equation}\label{mmbshpLqlam}
 (\lambda\rightarrow L(\cdot,\lambda))\in C(K,\;C^l(W,\;\mathfrak{B}(H^k))).
\end{equation}
By \cite[Lemma~7]{rz2012} and the arguments in the end of the proof of \cite[Lemma 8]{rz2012},
and in view of (\ref{boundLqmu}) and (\ref{mmbshpLqlam}), we get
 $(\lambda\rightarrow (\id+L(\cdot,\lambda))^{-1})\in C(K,\;C^l(\mathcal{W},\;\mathfrak{B}(H^k)))$.
Then (\ref{rprresHq}) and Lemma~\ref{lmschrhilb}(iv) imply the desired inclusion (\ref{inclRlambda}).
\qed

\begin{theorem}\label{contdeponqeigfunct}
Let $\lambda(q)$ be the least eigenvalue  of $\mathcal{H}_q\ (q\in\RR^n)$.
If $\beta\in C^\infty(F\times\RR^n)$ then $\lambda\in C^\infty(\RR^n)$ and
there exists a unique smooth section $e:\RR^n\mapsto L_2\times\RR^n$ such that
$e(\cdot,q)$ is a positive eigenfunction of $\mathcal{H}_q$
related to $\lambda(q)$ with $\Vert e(\cdot,q)\Vert_{L_2}=1$.
\end{theorem}

\noindent\textbf{Proof}. Assume without loss of generality that $\beta(x,q)<0$
(otherwise we can consider $\mathcal{H}_q-\mu$ instead of $\mathcal{H}_q$ with a suitable $\mu>0$).
Since $\lambda(q)$ is a simple eigenva\-lue of $\mathcal{H}_q$ for any $q\in\RR^n$,
there exists a unique positive eigenfunction $e(\cdot,q)$, related to it,
such that $\Vert e(\cdot,q)\Vert_{L_2}=1$.
Let $\lambda_0$ be the least eigenvalue of the operator $\mathcal{H}_{0}$
and $e_0$ be the eigenfunction related to
$\lambda_0$ and satisfying conditions mentioned above.
Let $\Gamma$ be a circle of small radius
in the complex plane $\CC$ not intersecting $\sigma(\mathcal{H}_0)$ and surrounding only~$\lambda_0$.
By~Lemma~\ref{lmcontres} (with $k=0$) one may restrict on the open
neighborhood $Q$ of $0$ in such a way that
$\Gamma\subset\CC\setminus\sigma(\mathcal{H}_q)$ for any $q\in Q$,
and the inclusion (\ref{inclRlambda}) is valid with $K=\Gamma$.
Hence, in view of $H^2\hookrightarrow L_2$, the Riesz projection
\begin{equation}\label{Rieszpr}
 P(q)=-\frac{1}{2\pi i}\oint_{\Gamma}R_\lambda(\mathcal{H}_q)\,d\lambda\quad (q\in Q)
\end{equation}
onto the invariant subspace of $\mathcal{H}_q$ corresponding to the
part of its spectrum lying inside of $\Gamma$ (\cite[Introduction, Sect.~4]{k72}) has the property for any $l\in\NN$:
\begin{equation}\label{mmbshpP}
 P\in C^l(Q, \mathfrak{B}(L_2)).
\end{equation}
In particular, one may restrict $Q$ in such a way  that
$\|P(q)-P(0)\|_{\mathfrak{B}(L_2)}\le 1/2$ for any $q\in Q$.
Then, taking into account that $P(q)$ are orthogonal projections
(since ${\mathcal H}_q$ are self-adjoint), we have $\dim({\rm Im}\,P(q))=\dim({\rm Im}\,P(0))$,
see \cite[Chapt.~III, Sect.~34]{Akh-Gl}.
This means that for $q\in Q$ the operator $\mathcal{H}_q$
has inside of $\Gamma$ only one simple eigenvalue $\tilde\lambda(q)$,
it is real because $\mathcal{H}_q$ is self-adjoint, ${\rm Im}\,P(q)$
is the eigenspace of $\mathcal{H}_q$ related to $\tilde\lambda(q)$,
and $\tilde\lambda(0)=\lambda_0$.
Denote $\tilde e(\cdot,q)=P(q)e_0$. We have for any $q\in Q$
\[
 \|\tilde e(\cdot,q)\|_{L_2}\ge\|P(0)e_0\|_{L_2}-\|(P(q)-P(0))e_0\|_{L_2}\ge{1}/{2}.
\]
Thus, $\tilde e(x,q)$ is an eigenvector of $\mathcal{H}_q$
related to $\tilde\lambda(q)$ such that $\tilde e(x,0)=e_0(x)$.
By~(\ref{mmbshpP}),
\begin{equation}\label{mmbshtile}
 (q\rightarrow\tilde e(\cdot,q))\in C^l(Q, L_2)
\end{equation}
for any $l\in\NN$. Then the equality
 $\tilde\lambda(q)\<\mathcal{H}_q^{-1}\tilde e(\cdot,q),\tilde e(\cdot,q)\>_{L_2}
 =\Vert\tilde e(\cdot,q)\Vert_{L_2}^2$
 and Lem\-ma~\ref{lmcontres} (with $k=0$ and $K=\{0\}$) imply that $\tilde\lambda\in C^\infty(Q,\RR)$.

Take an arbitrary $m\in\NN$ and set $j=[p/4+m/2]+1$.
Since the equality $\tilde e(\cdot,q)=(\tilde\lambda(q))^j\mathcal{H}_q^{-j}\tilde e(\cdot,q)$ is valid,
then in view of (\ref{mmbshtile}) and Lemma~\ref{lmcontres}, we can restrict $Q$ in such a way that
$(q\rightarrow\tilde e(\cdot,q))\in C^l(Q,\,H^{2j})$
for any $l\in\NN$.
On~the other hand, by Sobolev's Embedding Theorem, $H^{2j}\hookrightarrow C^m$, see \cite{aub}.
Thus, for any $m\in\NN$ there is an open neighborhood $Q$ of $0$ such that
 $(q\rightarrow\tilde e(\cdot,q))\in C^\infty(Q,\,C^m)$.
In particular, sin\-ce $e_0>0$, we can restrict $Q$ in such a way that
$e(x,q)={{\rm Re}(\tilde e(x,q))}/{\Vert{\rm Re}(\tilde e(\cdot,q))\Vert}_{L_2}>0$ for any $q\in Q$.
Clearly, $e(x,q)$ is an eigenfunction of $\mathcal{H}_q$, related to
the eigenvalue $\tilde\lambda_q$ and $\Vert e(\cdot,q)\Vert_{L_2}=1$.
 It remains only to show that it is possible to restrict $Q$ in such a way that
$\tilde\lambda(q)$ is the least eigenvalue of $\mathcal{H}_q$ for
any $q\in Q$, i.e., $\tilde\lambda(q)=\lambda(q)$.
Indeed, otherwise there is a sequence $q_\nu\in Q$ such that
$\lim_{\nu\rightarrow\infty}q_\nu=0$ and for any $\nu$ there
exists an eigenvalue $\tilde\lambda_\nu$ of $\mathcal{H}_{q_\nu}$ obeying
$\tilde\lambda_\nu<\tilde\lambda(q_\nu)$.
Since the operators $\mathcal{H}_{q_\nu}$ are positive definite and for some $\delta>0$
in the interval $(\lambda_0-\delta,\lambda_0+\delta)$ there is only
the eigenvalue $\tilde\lambda(q_\nu)$ of $\mathcal{H}_{q_\nu}$,
we get $\tilde\lambda_\nu\in[0,\lambda_0-\delta]$ for any $\nu$.
Let $\lambda_*\in[0,\lambda_0-\delta]$ be a concentration  point of the
sequence $\{\tilde\lambda_\nu\}_{\nu\in\NN}$.
Choosing a subsequence, we can assume that
$\lim_{\nu\rightarrow\infty}\tilde\lambda_\nu=\lambda_*$.
Surrounding $\lambda_*$ by a small enough circle $\Gamma$
such that $\Gamma\cap[\lambda_0,\infty)=\emptyset$, considering for
each $\nu$ the Riesz projection $\tilde P_\nu$, defined by the rhs
of (\ref{Rieszpr}) with $q=q_\nu$ and using the above arguments,
we get that $\lim_{\nu\rightarrow\infty}\Vert\tilde P_\nu-\tilde P(0)\Vert_{\mathfrak{B}(L_2)}=0$,
where $\tilde P(0)$ is defined by the rhs of (\ref{Rieszpr}) with $q=0$.
Since $q_\nu$ lies inside $\Gamma$ for a large eno\-ugh $\nu$, $\dim\mathrm{Im} P(0)>0$.
Hence, there is at least one eigenvalue of $\mathcal{H}_{0}$ inside of~$\Gamma$.
But~this is impossible, because $\lambda_0$ is the least eigenvalue of $\mathcal{H}_{0}$.
\qed

\smallskip
\textbf{3.3. Solution of the stationary equation}.
Consider the compact domain in $\RR\times F$
\begin{equation}\label{dfD}
 D:=\{(u,\,x)\in\RR\times F:\;u_-(x)\le u\le u_+(x)\},
\end{equation}
where $u_-,u_+\in C^\infty$ and $u_-\le u_+$.
Define sets $G^k=\mathrm{Int}(G)\cap C^k$ for $k>0$, where
$G\subset C$ is a bounded, closed and convex set given by
\begin{equation}\label{dfG}
 G:=\{u\in C:\ u_-(x)\le u(x)\le u_+(x) \ \ {\rm for\ all}\ \ x\in F\}.
\end{equation}

\begin{lemma}\label{lmcontsuper}
Let $m,l\in\NN$, $\Pi$ be an open domain in
$\RR^m\times F$ of the form $\Pi=\RR^{m-1}\times\mathrm{Int}(D)$,
and $\theta(\cdot,\cdot,\cdot):\Omega\times\RR^n\rightarrow\RR^l$ be a
continuous function. Then for any $q\in\RR^n$ the superposition
operator $\Theta(v,q)(x)=\theta(v(x),x,q)$ maps the set
\[
 Y=\{v\in C(F,\,\RR^m):\,(v(x),\,x)\in\Pi\ {\rm for\ all}\ x\in F\}
\]
into the set $C(F,\,\RR^l)$, and the inclusion $\Theta\in C(Y\times\RR^n,\,C(F,\,\RR^l))$ holds.
\end{lemma}

\noindent\textbf{Proof}. The first claim is obvious. Let us prove the second one.
Suppose that $v_0\in Y$. Consider the set $\Gamma(v_0)=\{(v_0(x),x)\}_{x\in F}$.
Take a relatively compact open set ${Q}\subset\RR^n$ such that $0\in {Q}$.
In order to construct a similar open neighborhood of the set $\Gamma(v_0)$,
observe that it is compact in~$\Pi$.
Then there is a finite open covering $\{U_j\}_{j=1}^k$ of $\Gamma(v_0)$ such that
$\overline{U_j}\subset\Pi\ (1\le j\le k)$ and each of $\overline{U_j}$ is compact.
Consider the open set $\Pi^\prime=\bigcup_{j=1}^k U_j$.
Then $\Gamma(v_0)\subset\Pi^\prime$,
$\overline{\Pi^\prime}=\bigcup_{j=1}^k\overline{U_j}\subset\Pi$
and the set $\overline{\Pi^\prime}$ is compact.
Consider the following open subset of $C(F,\,\RR^m)$:
\[
 Y^\prime=\{v\in C(F,\,\RR^m):\,(v(x),\,x)\in\Pi^\prime\ {\rm for\ all}\ x\in F\}.
\]
It is clear that $v_0\in Y^\prime$. Since $\theta(u,x,q)$ is uniformly continuous on the compact
$\overline{\Pi^\prime}\times\overline{{Q}^\prime}$,
for any $\epsilon>0$ there is $\delta>0$ such that
$|\theta(v_1,x,q)-\theta(v_2,x,0)|<\eps$ for all
$(v_1,\,x),\,(v_2,\,x)\in\overline{\Omega^\prime}$ and
$q\in\overline{{Q}^\prime}$, where $\|v_1-v_2\|<\delta$ and
$\|q\|<\delta$. Let us choose $\sigma\in(0,\delta)$ such that
 $B_\sigma(v_0)=\{v\in C(F,\,\RR^m):\,\Vert v-v_0\Vert_{C(F,\,\RR^m)}<\sigma\}\subset Y^\prime$.
If $v\in B_\sigma(u_0)$ and $\|q\|<\delta$ then
$\Vert\Theta(u,q)-\Theta(u_0,0)\Vert_{C(F,\,\RR^l)}<\epsilon$ holds.
\qed

\begin{lemma}\label{lmactsuper1}
Let $f\in C^\infty(D\times\RR^n)$. Then the superposition operator
\begin{equation}\label{E-superop}
 \Phi_f(u,q)=f(u(x),x,q)
\end{equation}
obeys $\Phi_f\in C^l(G^k\times\RR^n,\,C^k)$ for any integers $k,l\ge 0$.
\end{lemma}

\noindent\textbf{Proof}. This is divided into two steps.

Step 1. First, we shall reduce the operator
(\ref{E-superop}), acting from $G^k$ into $C^k$ to a
superposition operator, acting in spaces of continuous vector functions.
Take $(u,q)\in G^k\times\RR^n$ for some  $k\in\NN$.
Observe that differentials $d^ju(x)\,(1\le j\le k)$ can be conside\-red as functions
defined on $F$ with values in $\RR^{n_j}$ (e.g., $n_1=p$, $n_2=p(p+1)/2$ and so on).
We~have the following:
\begin{eqnarray*}
  && d_x\Phi_f(u,q) = \partial_x f(u(x),x,q)+\partial_u f(u(x),x,q)\,du(x),\\
 &&\! d_x^2\Phi_f(u,q) = \partial_x^2f(u(x),x,q) {+}\partial_u^2f(u(x),x,q)(du(x),du(x)) {+}\partial_u f(u(x),x,q) d^2u(x),
\end{eqnarray*}
and so on. Hence,
\begin{equation}\label{bigsuperpos}
\hskip-3mm
(d_x\Phi_f(u,q),\,d_x^2\Phi_f(u,q),\dots,\,d_x^k\Phi_f(u,q),\,\Phi_f(u,q))=
\psi(v(x),x,q),
\end{equation}
where $v(x)=(du(x),d^2u(x),\dots,d^ku(x),u(x))$ is the vector function on $F$,
and the function $\psi:\Omega\times\RR^n\rightarrow\RR^N$ is smooth.
Here $\Omega=\RR^{N-1}\times{\rm Int}\,D$ and $N=1+\sum_{j=1}^kn_j$.
It is enough to~show that the superposition operator $\Psi(v,q)=\psi(v(x),x,q)$ obeys
\begin{equation}\label{reducedmmmbship}
 \Psi\in C^l(X\times\RR^n,\;C(F,\,\RR^N))
\end{equation}
for any $l\in\NN$, where
 $X=\{v\in C(F,\,\RR^N):\,(v(x),\,x)\in\Omega\ {\rm for\ all}\ x\in F\}$.

Step 2. In aim to prove (\ref{reducedmmmbship}), take $v\in X$,
$q\in\RR^n$, $h\in C(F,\,\RR^N)$ and $\delta>0$ such that $u+th\in X$
for any $t\in[0,\delta]$. We have the following representation for
$t\in[0,\delta]$:
\begin{eqnarray*}
 t^{-1}\big(\Psi(v+th,q)-\Psi(v,q)\big) \eq t^{-1} h(x)\int_0^t\partial_v\psi(v(x)+\tau h(x),x,q)\dtau.
\end{eqnarray*}
Hence,
\begin{eqnarray*}
 &&\Vert t^{-1}(\Psi(v+th,q)-\Psi(v,q))-\partial_v\psi(u(x),x,q)\,h(x)\Vert_{C(F,\,\RR^N)}\le\\
 &&\le
 \Vert h\Vert_{C(F,\,\RR^N)}\sup\nolimits_{\,\tau\in[0,t]}
 \Vert\partial_v\psi(v +\tau h,\cdot,q)-\partial_v\psi(v,\cdot,q)\Vert_{C(F,\,\mathfrak{B}(\RR^N))}.
\end{eqnarray*}
Since for any fixed $q\in\RR^n$ the function $g(x,\tau)=\partial_v\psi(v(x)+\tau h(x),x,q)$ is uniformly
continuous in $F\times[0,\delta]$, the last estimate implies
$$
 \lim\nolimits_{\,t\downarrow 0} \Vert\,t^{-1}g(\Psi(v+th,q)-\Psi(v,q))-\partial_v\psi(v(x),x,q)\,h(x)\Vert_{\,C(F,\,\RR^N)}=0.
$$
Hence, for any $q\in\RR^n$ the operator $\;\Psi(\cdot,q):\,X\rightarrow C(F,\RR^N)$ has at any $u\in X$
the Gateaux partial differential $\partial_v\Psi(v,q)$ of the form
$$
 \partial_v\Psi(v,q)\,h=\partial_v\psi(v(x),x,q)\,h(x),\quad h\in C(F,\RR^N).
$$
We see that $\partial_v\Psi(v,q)\in\mathfrak{B}(C(F,\,\RR^N))$ and
\begin{eqnarray*}
 &&\Vert\partial_v\Psi(v_1,q_1)-\partial_v\Psi(v_2,q_2)\Vert_{\mathfrak{B}(C(F,\,\RR^N))}\\
 &&\quad\le\Vert\partial_v\psi(v_1,\cdot,q_1)-\partial_v\psi(v_2,\cdot,q_2)\Vert_{C(F,\,\mathfrak{B}(\RR^N))}.
\end{eqnarray*}
By Lemma~\ref{lmcontsuper} applied to the superposition operator
\[
 \Psi_1(v,q)=\partial_v\psi(v(x),x,q):\;C(F,\,\RR^N)\rightarrow C(F,\,\mathfrak{B}(\RR^N)),
\]
the partial differential $\partial_v\Psi(v,q)$ is continuous in the sense that
 $\partial_v\Psi$ belongs to $C(C(F,\,\RR^N)\times\RR^n,\;\mathfrak{B}(C(F,\,\RR^N)))$.
Hence, $\partial_v\Psi(v,q)$ is the Freshet's partial differential.
Similarly, one may show that at any $(v,q)\in X\times\RR^n$
there is the Gateaux partial differential
$$
 \partial_q\Psi(v,q)\,s=\lim\nolimits_{\,t\downarrow 0}\,t^{-1}(\Psi(v,q+ts)-\Psi(v,q))
 =\partial_q\psi(v(x),x,q)\,s\,,
$$
where $s\in\RR^n$, and the limit is taken with respect to the $C(F,\,\RR^N)$-norm.
Furthermore, this differential is continuous:
$
 \partial_q\Psi\in C(C(F,\,\RR^N)\times\RR^n,\;\mathfrak{B}(\RR^n,\,C(F,\,\RR^N))).
$
Again, this fact implies that $\partial_q\Psi(v,q)$ is the Freshet's partial differential.
Thus, we have proved that the superposition operator $\Psi(v,q)(x)=\psi(v(x),x,q)$ belongs to class
$C^1(X\times\RR^n,\,C(F,\,\RR^N))$.
Applying similar arguments to superposition operators
\[
 \Psi_1(v,q)=\partial_v\psi(v(x),x,q),\quad
 \Psi_2(v,q)=\partial_q\psi(v(x),x,q),
\]
one may show that $\Psi\in C^2(X\times\RR^n,\,C(F,\,\RR^N))$.
Further, we can prove by induction that (\ref{reducedmmmbship}) holds for any $l\in\NN$.
\qed

\smallskip

 For $\alpha\in (0,1)$ and integer $k\ge 0$ consider the Banach space
$C^{k,\,\alpha}$ of such functions $u\in C^k$, for which all
the partial derivatives of order $k$ belong to H\"older class
$C^{0,\,\alpha}$. The norm in this space is defined as follows:
\begin{equation}\label{dfCkalphanrm}
 \Vert u\Vert_{C^{k,\,\alpha}}=\max\big\{\max_{|\beta|\le k}\Vert D^\beta u\Vert_{C^0},\;
 \max_{|\beta|=k}\;\sup_{x,y\in F,\;x\neq y} {|D^\beta u(x)-D^\beta u(y)|}\,{d(x,y)^{-\alpha}}\big\}.
\end{equation}

\begin{lemma}\label{lmsmoothellipt}
Let $\Delta_q$ be the elliptic operator in (\ref{ellop}).
Then for any $l\in\NN$ the mapping $\mathfrak{D}:(u,q)\mapsto-\Delta_q\,u$ belongs to class
$C^l(C^{k+2,\,\alpha}\times\RR^n,\;C^{k,\,\alpha})$.
\end{lemma}

\noindent\textbf{Proof}.
Fix a finite atlas $\{(U_a,x_a)\}_{1\le a\le A}$ of $F$.
Taking $u\in C^{k+2,\,\alpha}$, $q\in\RR^n$ and $s\in\RR^n$, we obtain in a local chart, see (\ref{ellop}):
\begin{eqnarray*}
 && 1/t(\mathfrak{D}(u,q+ts)-\mathfrak{D}(u,q))-\mathfrak{D}_1(u,q)\,s\\
 && = -\big(1/t\int_0^t\partial_q g^{ij}(x,q+\tau s)\,s\dtau
 -\partial_q g^{ij}(x,q)\,s\big)\partial^2_{ij}\,u\\
 && - \big(1/t\int_0^t\partial_q b^i(x,q+\tau s)\,s\dtau
 -\partial_q b^i(x,q)\,s\big)\,\partial_i u,
\end{eqnarray*}
where
 $\mathfrak{D}_1(\cdot,q)\,s = -\partial_q g^{ij}(x,q)\,s\partial^2_{ij}
 -\partial_q b^i(x,q)\,s\partial_i$.
In view of (\ref{dfCkalphanrm}),
\begin{eqnarray*}
 &&\Vert 1/t(\mathfrak{D}(u,q+ts)-\mathfrak{D}(u,q))-\mathfrak{D}_1(u,q)\,s\,\Vert_{C^{k,\,\alpha}(U_a)}\\
 &&\le C_a\max_{i,j\in\{1,2,\dots,p\}}\max_{\tau\in[0,t],\,x\in U_a}
   \big(
   |\partial_qg^{ij}(x,q+\tau s)-\partial_q g^{ij}(x,q)|\\
 &&+\|\partial_q b^i(x,q+\tau s)-\partial_q b^i(x,q)\|_{C^{k+1}(U_\alpha,{\cal B}(\RR^n))}\big)|s|\cdot\Vert u\Vert_{C^{k+2,\,\alpha}(U_a)}
\end{eqnarray*}
holds for some $C_a>0$ that does not depend on $u$.
Replacing $C_a$ by $C=\max_{1\le a\le A} C_a$, we find that
$\mathfrak{D}:\,C^{k+2,\,\alpha}\times\RR^n\rightarrow C^{k,\,\alpha}$
has the partial Gateaux differential $\partial_q\mathfrak{D}(u,q)$ at each point $(u,q)$,
and it is equal to $\mathfrak{D}_1(u,q)$.
Similarly to the proof of Lemma~\ref{lmsmoothellipt1}, we obtain that $\mathfrak{D}$ has
continuous mixed partial differentials by $q$ and $x$ (of any order)
at any point $(u,q)\in C^{k+2,\,\alpha}\times\RR^n$.
\qed

\begin{theorem}\label{prsmoothparstat}
Let $f\in C^\infty(D\times\RR^n)$ and $u_*(x)\in\mathrm{Int}\,G$ be a smooth solution of (\ref{nonlinstatpar1}) with $q=0$
such that $\lambda=0$ is not an eigenvalue of the
operator ${\mathcal H}=-\Delta-\partial_u f(u_*(x),x,0)$ on $L_2$ with domain in $H^2$.
Then for any integers $k\ge 0$ and $l\ge 1$ and $\alpha\in(0,1)$ there are open neighborhoods
$U_*\subseteq C^{k+2,\alpha}$ of $u_*$ and $V_0\subseteq\RR^n$ of $\,0$
such for any $q\in V_0$ there exists in $U_*$ a unique solution
$u(x,q)$ of (\ref{nonlinstatpar1}), in particular, $u_*(x)=u(x,0)$,
such that the function $q\rightarrow u(\cdot,q)$ belongs to class $C^l(V_0,U_*)$.
\end{theorem}

\noindent\textbf{Proof}. By Lemma~\ref{lmactsuper1}, for any  $q\in\RR^n$ and
integers $k,l\ge0$ the operator (\ref{E-superop}) maps the
set $G_{k+2}=\mathrm{Int}\,G\cap C^{k+2}$ into $C^{k+2}$ and
 $\Phi_f\in C^l(G_{k+2}\times\RR^n,\,C^{k+2})$.
Since $C^{k+2,\alpha}$ and $C^{k+2}$ are continuously embedded
into $C^{k+2}$ and $C^{k,\alpha}$, respectively, we get
 $\Phi_f\in C^l(G_{k+2,\alpha}\times\RR^n,\,C^{k,\alpha})$,
where $G_{k+2,\alpha}=\mathrm{Int}\,G\cap C^{k+2,\alpha}$.
Consider operators
 $Y(u,q)=\Delta_q\,u+\Phi_f(u,q)\ (q\in\RR^n)$
defined on $G_{k+2,\,\alpha}$.
By Lemma~\ref{lmsmoothellipt},
$
 Y\in C^l(G_{k+2,\alpha}\times\RR^n,C^{k,\alpha}).
$

Let $\widetilde{\mathcal H}$ be ${\mathcal H}$ restricted
on $C^{k,\,\alpha}$ with the domain $C^{k+2,\,\alpha}$.
Set~$\beta=\partial_u f(u_*(x),x,0)$.
By (\ref{dfCkalphanrm}), there is $C>0$
such that $\Vert\Delta u\Vert_{C^{k,\,\alpha}}\le C\Vert u\Vert_{C^{k+2,\,\alpha}}$ for any $u\in C^{k,\,\alpha}$;
hence, $\widetilde{\mathcal H}-\mu\in\mathfrak{B}(C^{k+2,\,\alpha},\,C^{k,\,\alpha})$.
Let $\mu<-\max_{\,x\in F}\,\beta(x)$.
By~\cite[Theorem~4.18]{aub}, ${\cal H}(u)-\mu\,u=f(x)$
has a unique solution $u\in C^{k+2,\,\alpha}$ for any $f\in C^{k,\,\alpha}$,
i.e., $\widetilde{\mathcal H}-\mu$ maps $C^{k+2,\,\alpha}$ injectively onto $C^{k,\,\alpha}$.
As in the proof of Lemma~\ref{lmschrhilb}(i), using Banach's Closed
Graph Theorem and compactness of the embedding $C^{k+2,\,\alpha}\hookrightarrow C^{k,\,\alpha}$,
we prove that $\widetilde{\mathcal H}-\mu$ is continuously invertible
and $(\widetilde{\mathcal H}-\mu)^{-1}:\;C^{k,\,\alpha}\rightarrow C^{k,\,\alpha}$ is compact.

By the above, the spectrum of $(\widetilde{\mathcal H}-\mu)^{-1}$ consists
of the point $\nu=0$ and a countable number of non-zero eigenvalues $\nu_n$
of finite multiplicity, which can accumulate only at the point $\nu=0$.
Hence the spectrum of $\widetilde{\mathcal H}$ is discrete.
By Lemma~\ref{lmschrhilb}(i-iii), $\sigma({\mathcal H})$ is discrete and $\mu\notin\sigma({\mathcal H})$.
By the same arguments as in the proof of Lemma~\ref{lmschrhilb}(iii),
we find that the spectrum of $\widetilde{\mathcal H}$ is discrete
and coincides with $\sigma({\mathcal H})$.

Let $\lambda\notin\sigma({\mathcal H})$.
Then $\widetilde{\mathcal H}-\lambda$
maps injectively $C^{k+2,\,\alpha}$ into $C^{k,\,\alpha}$ and it is bounded.
By the Elliptic Regu\-larity Theorem, this operator is surjective.
By the Banach's Closed Graph Theorem, $\widetilde{\mathcal H}-\lambda$ acts from
$C^{k+2,\,\alpha}$ into $C^{k,\,\alpha}$ and it is continuously invertible.
Since $0\notin\sigma({\mathcal H})$, the partial differential is continuously invertible
\[
 \partial_u Y(u_*,0) =-{\cal H}
 \in
 \mathfrak{B}(C^{k+2,\,\alpha},\,C^{k,\,\alpha}).
\]
These facts and the Implicit Function Theorem, see \cite{aub}, complete the proof.
\qed


\end{document}